\documentclass[12pt,reqno]{amsart}
\usepackage{amsmath}
\usepackage{amssymb}
\usepackage[left=1in,top=1.4in,right=1in,bottom=1.2in]{geometry}
\usepackage[usenames]{color}
\usepackage[dvipsnames]{xcolor}
\usepackage{enumerate}
\usepackage[normalem]{ulem}
\usepackage{soul}
\usepackage{graphicx, caption, subcaption}

\usepackage{tikz}
\usepackage{fp}  
\usetikzlibrary{positioning,chains,fit,shapes,calc}
\usepackage{subcaption}

\newcommand\Item[1][]{%
  \ifx\relax#1\relax  \item \else \item[#1] \fi
  \abovedisplayskip=0pt\abovedisplayshortskip=0pt~\vspace*{-\baselineskip}}

\newtheorem{theorem}{Theorem}[section]
\newtheorem{lemma}[theorem]{Lemma}

\newtheorem{conjecture}[theorem]{Conjecture}

\def\C{k}

\def\nn{\nonumber}
 \def\b{\beta}  \def\D{\Delta}
    
  \def\k{\kappa}

\newcommand{\rbrac}[1]{\left(#1\right)}
\newcommand{\sbrac}[1]{\left[ #1\right]}
\newcommand{\cbrac}[1]{\left\{ #1\right\}}

\def\E{\mathbb{E}}

\def\Var{\mbox{{\bf Var}}}
\def\P{\mathbb{P}}

\def\codeg{\text{codeg}}

\newcommand{\of}[1]{\left( #1 \right) }

\newcommand{\sqbs}[1]{\left[ #1 \right]}

\newcommand{\Mean}[1]{\E\sqbs{#1}}

\allowdisplaybreaks[1]
\newcommand{\ignore}[1]{}

\newcommand{\beq}[1]{\begin{equation}\label{#1}}
\newcommand{\eeq}{\end{equation}}

\newcommand{\tO}[1]{\tilde{O}\rbrac{#1}}

\newcommand{\mc}[1]{\mathcal{#1}}

\def\hati{\hat{i}}
\def\hatt{\hat{t}}

\title[Closing the Random Graph Gap in Tuza's Conjecture]{Closing the Random Graph Gap in Tuza's Conjecture Through the Online Triangle Packing Process}

\author{Patrick Bennett}
\address{Department of Mathematics, Western Michigan University, Kalamazoo, MI, USA}
\thanks{The first author was supported in part by Simons Foundation Grant \#426894.}
\email{\tt patrick.bennett@wmich.edu}

\author{Ryan Cushman}\thanks{}
\address{Department of Mathematics,
Western Michigan University, Kalamazoo, MI, USA} 
\email{ryan.cushman@wmich.edu}

\author{Andrzej Dudek}
\address{Department of Mathematics, Western Michigan University, Kalamazoo, MI, USA}
\thanks{The third author was supported in part by Simons Foundation Grant \#522400.}
\email{\tt andrzej.dudek@wmich.edu}

\begin{document}

\begin{abstract}
A long-standing conjecture of Zsolt Tuza asserts that the triangle covering number $\tau(G)$ is at most twice the triangle packing number $\nu(G)$,  where the \emph{triangle packing number} $\nu(G)$ is the maximum size of a set of edge-disjoint triangles in $G$ and the \emph{triangle covering number} $\tau(G)$ is the minimal size of a set of edges intersecting all triangles. In this paper, we prove that Tuza's conjecture holds in the Erd\H{o}s-R\'enyi random graph $G(n,m)$ for all range of $m$, closing the ``gap'' in what was previously known. (Recently, this result was also independently proved by Jeff Kahn and Jinyoung Park.)
We employ a random greedy process called the \emph{online triangle packing process} to produce a triangle packing in $G(n,m)$ and analyze this process by using the differential equations method.
\end{abstract}

\maketitle

\section{Introduction}

Random processes that seek to control the presence of triangles in a graph have been extensively and fruitfully applied to a number of important problems. The best known of these is the \emph{triangle-free process} introduced by Bollob\'as and Erd\H{o}s (see \cite{BR}), which maintains a triangle-free subgraph $G_T(i) \subset G(n,i)$. Here $G(n,m)$ is the Erd\H{o}s-R\'enyi random graph that assigns equal probability to all graphs on a fixed set $V$ of $n$ vertices with exactly~$m=m(n)$ edges. At each step in the triangle-free process, an edge is revealed and added to $G_T(i)$ only if it does not create a triangle in $G_T(i)$. Famously used to study the Ramsey numbers $R(3,t)$, the triangle-free process has continued to yield results in this area; indeed, this can be seen from the recent results of Bohman and Keevash \cite{BK2} and independently of Fiz Pontiveros, Griffiths and Morris \cite{FGM}, which gives the best-known lower bound of \newline $R(3, t) \ge \rbrac{1/4-o(1)}t^2 / \log t$. 

Another important process, called \emph{random triangle removal} and also introduced by Bollob\'as and Erd\H{o}s, creates a triangle-free graph from a complete graph instead of an empty graph (see \cite{B1, B2}). Here we start with $G_R(0)=K_n$ and remove all three edges from a triangle chosen uniformly at random from the triangles in $G_R(i)$ until every triangle is removed. Note that the edges of the removed triangles form a triangle packing in $K_n$. Although originally motivated by the study of $R(3,t)$, this process has not yet resulted in any good bound. Bohman, Frieze and Lubetzky \cite{BFL} provided the best-known upper and lower bounds of $n^{3/2 + o(1)}$ on the number of edges with high probability remaining at the end of this process. (We say an event dependent on $n$ occurs \emph{with high probability} (abbreviated w.h.p.) if the probability of that event tends to one as $n$ tends to infinity.) 

A third process introduced by Bollob\'as and Erd\H{o}s is the \emph{reverse triangle-free process}. Similar to the previous process, this one starts with $G_{RT}(0) = K_n$ but instead removes \emph{one} edge that is in a triangle from $G_{RT}(i)$. As before, the process terminates when the graph is triangle-free. The number of edges in the final graph is w.h.p. $(1+o(1))\sqrt{\pi} n^{3/2}/4$ due to Erd\H{o}s, Suen and Winkler~\cite{ESW} and the final number of edges is concentrated about its expectation due to Makai~\cite{Makai} and independently Warnke~\cite{LutzRev}.

In this paper, we find a triangle packing in $G(n,m)$ by using a process inspired by this rich history: the online triangle packing process. Our application of this process will allow us to improve the best-known bound on the triangle packing number for random graphs. We define the \emph{triangle packing number} of a graph $G$, denoted $\nu(G)$, as the maximal size of a set of edge-disjoint triangles. Here, we start with an empty packing $M(0)$ in $G(n,0)$ and reveal one edge at a time. If this edge forms a triangle that is edge-disjoint from $M(i)$, then many triangles might be created as well. Thus we choose one of the triangles uniformly at random from the set of created triangles and add its edges to $M(i+1)$. Hence the \emph{unmatched graph} $U(i) = G(n,i) - M(i)$ has no triangles by induction on $i$. (We will find it convenient in this paper to identify a graph $H$ with its edge set $E(H)$.) Further, as the name implies, the triangles of $M(i)$ form a triangle packing. 

Our analysis of the online triangle packing process will be done in a similar manner as was productively employed to study the other triangle-controlling processes described above: \emph{dynamic concentration} (also called the \emph{differential equation method}, see \cite{nick2, BD}). In this method, a system of random variables are tracked using martingale concentration inequalities to show that these variables w.h.p. stay close to what we expect them to be. 

The version of the online triangle packing process we use here is a refinement of the one used in by Bennett, Dudek, and Zerbib in~\cite{BDZ}. Recall that the creation of a triangle in $U(i)$ potentially coincides with the creation of many triangles. More specifically, it coincides with the creation of a copy of the tripartite graph $K_{1,1,s}$ for some $s \ge 1$. In~\cite{BDZ}, instead of choosing a triangle uniformly at random from the set of created triangles, the edges of $K_{1,1,s}$ were moved to $M(i+1)$ for  maximal $s\ge 1$. Then $U$ is triangle-free and a triangle packing can be obtained by choosing one triangle from each copy of $K_{1,1,s}$ in $M$. This modification was done to aid in technical details of the analysis, but also resulted in more edges being moved to the matched set at each step in the greedy algorithm than was necessary to obtain a triangle packing. In contrast, our form of the process moves only what is necessary to have $U$ be triangle-free and $M$ be a triangle packing. In addition, we were able to simplify the troublesome technical details, resulting in a more streamlined analysis.  

It is not surprising, then, that we were able to improve the lower bound on $\nu(G(n,m))$. Note that we are only concerned with the sparse case of $G(n,\C m)$ for $\C < (\log n)^2$ since Frankl and R\"odl~\cite{FR} gave a bound on $\nu(G)$ that is optimal in order for $\C \ge (\log n)^2$. (This was slightly improved by Pippenger by decreasing the bound on $\C$, see \cite{AY}.) In \cite{BDZ} the differential equation $z' = 2e^{-z^2}- 4z^2$ was used to model the unmatched degree heuristically given as $z(t)n^{1/2}$ for $t\ge 0$. Several useful properties of $z$ were determined, including that $0 \le z(t) \le .5932$.  If we define $L_\nu(\C)$ as
$$
  L_\nu(\C):=  \frac13 \sbrac{\C - \frac{z(\C)}{2} - 2\int_{0}^{\C} \sbrac{z(t)^2-1+e^{-z(t)^2}}\, dt},
$$
we may state the main result using this previous iteration of the online triangle packing process as follows: for arbitrarily small $\varepsilon > 0$, if $n^{-(1/20)+\varepsilon} <  \C \le \frac{1}{1000}\log \log n$ then w.h.p $$\nu(G(n,\C n^{3/2}))\ge (1+o(1))L_\nu(\C)n^{3/2}.$$ We contrast this with our main result. As before, we rely on a differential equation in $y(t)$ for $t\ge 0$ that models the unmatched degree $y(t)n^{1/2}$, given by $y'=6e^{-y^2}-4$. The function $y(t)$ is discussed more fully in Section~\ref{proof:preliminaries}, where it will be shown that $0 \le y(t)\le .6368$. With this in mind, our main result is the following.

\begin{theorem}\label{thm:main}
Let $G=G(n,m)$ be a random graph of order $n$ and size $m=\C n^{3/2}$, where $0 \le \C < \frac{1}{10000} \log \log n$.
Then, w.h.p.
\[
\nu(G) \ge (1+o(1))\frac{1}{3}\left(\C - \frac{y(\C)}{2} \right)  n^{3/2}.
\]
\end{theorem}
The proof of Theorem~\ref{thm:main} can be found in Section~\ref{sec:main_proof}. The improvement of the previous result can be seen in the removal of the integral term in the definition of $L_\nu$, which was necessary for accounting for the edges placed into $M$ that were not used in the triangle packing. Recall that these unnecessary edges from  the previous process were a result of moving copies of $K_{1,1,s}$ for $s\ge 1$ into $M$, but only using the three edges of a triangle for each such copy and disregarding the rest. In our current result this adjustment is unnecessary, resulting in a larger lower bound. In fact, this bound is good enough to conclusively finish the proof of a long-standing conjecture of Zsolt Tuza in the case of random graphs. 

Tuza's conjecture relates the triangle packing number and the triangle covering number. The \emph{triangle covering number} $\tau(G)$ is the minimal size of a set of edges intersecting all triangles. It is easy to see that $\nu(G) \le \tau(G) \le 3 \nu(G)$ for any graph $G$. Tuza, however, conjectured that this trivial upper bound could be lowered.

\begin{conjecture}[Tuza~\cite{tuza}]\label{tuzaconj} 
 For every graph $G$, $\tau(G) \le 2 \nu(G)$.
\end{conjecture}

The conjecture is tight for the complete graphs of order 4 and 5. The best-known upper bound is $\tau(G) \le \frac{66}{23} \nu(G)$ from Haxell~\cite{haxell}. A recent development is due to Baron and Kahn~\cite{BK}, who showed that, in general, the multiplicative constant 2 in Tuza's conjecture cannot be improved. They demonstrated   that for any $\alpha>0$ there are arbitrarily large graphs $G$ of positive density satisfying $\tau(G) > (1-o(1))|G|/2$ and $\nu(G) < (1+\alpha)|G|/4$, disproving a conjecture of Yuster~\cite{yuster}. See \cite{krivelevich, HR, AZ} for related results.

We will modify the analysis of the triangle-free process from Bohman~\cite{bohman} in order to obtain an upper bound on $\tau(G(n,m))$ and pair this with our lower bound on $\nu(G(n,m))$. Note that we are concerned here with a small range of $m$ for which Tuza's conjecture is still open. Using the bound for $\nu(G(n,m))$ discussed above, Bennett, Dudek and Zerbib \cite{BDZ} proved the conjecture for $G(n,m)$, with the exception of a small range of $m$.

\begin{theorem}[Bennett, Dudek, and Zerbib~\cite{BDZ}]\label{thm:tuza-old}
There  exist absolute constants $0<c_1<c_2$ such that if $m \le c_1n^{3/2}$ or $m\ge c_2n^{3/2}$, then w.h.p.~Tuza's conjecture holds  for $G=G(n,m)$.
\end{theorem}

The proof of Theorem~\ref{thm:tuza-old} gives that one can take $c_1:=0.2403$ and $c_2 := 2.1243$. The existence of the constant $c_1$ was recently also proved by Basit and Galvin~\cite{BG}. This ``gap'' in the values of $m$ for which Tuza's conjecture holds is due to the ``wastefulness'' of the previous iteration of the online triangle packing process. Our current iteration of the online triangle packing process and its subsequent bound on $\nu(G(n,m))$ will be enough for us to close this gap.

\begin{theorem}\label{thm:tuza}
Tuza's conjecture holds w.h.p. for~$G(n,m)$ for any range of $m$.
\end{theorem}
Theorem~\ref{thm:tuza} will be proved in Section~\ref{sec:tuza}. We recently learned that Theorem~\ref{thm:tuza} was independently proved by Kahn and Park~\cite{KP} using a very different approach.


\section{Finding a triangle packing through the random process}\label{sec:main_proof}

\subsection{Outline of the algorithm}

Our process reveals one edge of $G(n, m)$ at each step. So for step $i$ we have $G(n,i)$, whose edges we partition into subgraphs:  the {\em matched graph} $M(i)$ and the {\em unmatched graph} $U(i)$. We will maintain the property that $U(i)$ is triangle-free, and $M(i)$ is the union of disjoint triangles. Let $e_i$ be the edge we add at step $i$. 
If $U(i) \cup \{e_i\}$ is triangle-free then we let $U(i+1):= U(i) \cup \{e_i\}$ and $M(i+1):=M(i)$. Otherwise $U(i) \cup \{e_i\}$ has at least one triangle (any such triangle must use $e_i$ since $U(i)$ is triangle-free), and we choose one such triangle uniformly at random.  If from among the set of such triangles, we choose say $T$ (a set of three edges), then we set  $M(i+1):=M(i) \cup T$ and  $U(i+1)=U(i) \setminus T$. We remark that the edge $e_i = \lbrace u, v \rbrace$ creates a triangle precisely when the codegree of $u$ and $v$ in $U(i)$ is positive. Recall that the {\em codegree} of two vertices $u$ and $v$ in a graph $H$, written $\codeg_H(u,v)$, is the number vertices $w$ such that both $uw$ and $vw$ are edges of $H$.

We are also concerned with the degree of each vertex in both $M(i)$ and $U(i)$. Fix a vertex $v$. Then we write  $d_U(v, i)=\deg_{U(i)}(v)$ and $d_M(v, i)=\deg_{M(i)}(v)$ to represent the unmatched and matched degree at step~$i$, respectively. We will also write $d_G(v, i)= \deg_{G(n,i)}(v) = d_U(v, i)+d_M(v, i)$. For convenience, we will sometimes suppress ``$i$'' in this notation when it is clear from context. 

Now define the scaled time parameter as
\[
t = t(i) := \frac{i}{n^{3/2}}
\]
for $0\le i \le \frac{1}{10000} n^{3/2} \log \log n$. At each step $i+1$ we choose a random edge without replacement. Thus, the probability of choosing any particular edge that has not been chosen yet is
\[
 \frac{1}{\binom{n}{2} -i} = \frac {2}{n^2} (1+\tilde{O}(n^{-1/2})),
\]
where $a(n) \in \tilde{O}(b(n))$ if there exists $k\ge 0$ such that $a(n) \in O(b(n)\log^k b(n))$.

Next we describe some important heuristics which will be formally justified later. The first of these is that at each step~$i$ (excluding steps near the start), we have $$d_U(v) + d_M(v) =   deg_{G(n,i)}(v) = \frac{2i}{n} (1+o(1)) = 2tn^{1/2} (1+o(1))$$ for sufficiently large $m$. This is due to the concentration of vertex degrees in $G(n, m)$. Assuming heuristically that $d_U(v) \approx y(t)n^{1/2}$ and that the codegrees in $U(i)$ are distributed Poisson with expectation $n(yn^{-1/2})^2 =y^2$, we may say that $d_M(v) \approx (2t-y(t))n^{1/2}$ and that the number of unmatched edges is about $\frac12n^{3/2}y$.

With this framework, we calculate the one-step change in the number of unmatched edges. If $e_i = \lbrace u, v \rbrace$, we gain one unmatched edge when $\codeg_U(u,v) = 0$, which occurs with probability $e^{-y^2}$. If the unmatched codegree is positive, then we lose two unmatched edges; they move into the matched graph along with $e_i$. Thus, approximating the one-step change as a derivative, we arrive at the following:

\[
\Delta\rbrac{\frac12 y(t) n^{3/2}} \approx \rbrac{\frac12 y'(t) n^{3/2}}\Delta t = \frac 12 y' \approx 1\cdot e^{-y^2} - 2 \cdot (1-e^{-y^2}).
\]
Here we use the fact that $\Delta t = n^{-3/2}$. Hence, we get the differential equation $y' = 6e^{-y^2}-4$. (We discuss this differential equation further in Section~\ref{proof:preliminaries}.)

To conclude our outline of the algorithm, we recall that the number of matched edges is $kn^{3/2}-\frac{y(k)}{2}n^{3/2}$ after $kn^{3/2}$ edges have been revealed. Thus the number of edge-disjoint triangles we obtain at the end of our process should be 
$$
\left(\frac13 k - \frac 16 y(k)\right)n^{3/2}.
$$
We show that with high probability this is very close to the actual situation.

\subsection{Preliminaries}\label{proof:preliminaries}
Let $y=y(t)$ for $t\ge 0$ be such that  the following autonomous differential equation holds:
\begin{equation*}\label{eq:zdiff}
y'=6e^{-y^2}-4.
\end{equation*}
 Assume that $y(0)=0$.  Then $y$ is {an increasing function of} $t$ and $y$ approaches  the unique positive root of the equation $6e^{-x^2}-4=0$  (as $t$ goes to infinity), which is  $\zeta = \sqrt{\log\rbrac{\frac32}} \approx0.6367$. Hence, $0\le y\le \zeta$. This also implies that $y'(t)\ge 0$.

Furthermore, note that
\begin{equation}\label{eq:z2}
y'' = -12e^{-y^2}yy'\le 0
\end{equation}
and consequently $0\le y'\le y'(0)=2$.

For integers $b,c\ge 0$ let us define the following random variables for every step $i\ge 0$:
\begin{itemize}
\item $Q_{b,c}(u,v)=Q_{b, c}(u, v,i)$ is the set of vertices $w$ such that  $\codeg_U(w, u)=b$ and $\codeg_U(w, v)=c$.
\item $R_c(v)=R_c(v,i)$ is the set of vertices $u$ such that $\codeg_U(u, v)=c$.
\item $S_{c}(u,v)=S_{c}(u, v,i)$ is the set of vertices $w\in N_U(v)$ such that the unmatched codegree of $w$ and $u$, excluding $v$, is $c$.
\end{itemize}

We will sometimes write the name of a set when we mean the cardinality of that set.
Next we wish to define deterministic counterparts to these random variables. Here we use our heuristic that the unmatched graph is almost regular with degree $yn^{1/2}$ and that the codegrees are almost independent Poisson variables with expectation $y^2$. So for example for $Q_{b,c}$, we have
$$
\sum_{w}\Pr\left(\codeg_U(w,u) = b \text{ and }\codeg_U(w,v) = c\right) \approx \left(\frac{e^{-y^2}y^{2b}}{b!}\right)\left(\frac{e^{-y^2}y^{2c}}{c!}\right)n =\frac{e^{-2y^2}y^{2b+2c}}{b!c!}n.
$$
Reasoning in a similar way, we may define the following functions: 
\[
q_{b, c}= q_{b, c}(t):=\frac{e^{-2y^2}y^{2b+2c}}{b!c!},\quad\quad\quad r_c= r_c(t):=\frac{e^{-y^2}y^{2c}}{c!}, \quad \quad \quad s_c = s_c(t) := \frac{e^{-y^2}y^{2c+1}}{c!}.
\]
(Here we scale by an appropriate power of $n$.)

Observe that when $b=c=0$ we have $q_{0,0}=e^{-2y^2}$, $r_0 = e^{-y^2}$, and $s_0=e^{-y^2}y$. Moreover, 
since for any $k \ge 0$ and $0 \le x \le 1$, we get $e^{-x^2}x^k \le 1$, we obtain
\begin{equation}
q_{b, c}\le\frac{1}{b!c!}, \quad \quad \quad r_c \le \frac{1}{c!}, \quad \quad \quad  s_c\le  \frac{1}{c!}. \label{eq:cpqt-bnd}
\end{equation}

 Define an ``error function" 
\[
f(t):= \exp\cbrac{\frac{1000\log n}{\log \log n} \cdot t} n^{-1/5}
\]
and observe that for $0 \le t \le \frac{1}{10000} \log \log n$ we have $n^{-1/5} \le f(t) \le n^{-1/10}$.

Now we define the ``good event'' at step $i$. For a given step $i$, let $\mc{E}_i$ be the event such that  in $G=G(n,i)$ we have:
\begin{enumerate}[(i)] 
\item \label{nobigcod} \emph{No huge codegree:}  for all $u, v\in V$ we have 
\[
\codeg_{G}(u, v) \le \frac{3\log n}{\log \log n}=: c_{max}.
\] 


\item \label{dynconc} \emph{Dynamic concentration:}  for every $j \le i$,
\begin{itemize}
\item $\displaystyle  d_{G}(v,j) \in \of{2t \pm n^{-1/4} \log^2 n}n^{1/2}$, 
\item $\displaystyle  d_U(v,j) \in  \of{y\pm f  }n^{1/2}$,
\item $\displaystyle  |Q_{b, c}(u,v,j)| \in \of{q_{b, c}\pm  f } n$,
\item $\displaystyle  |R_c(v,j)| \in  \of{r_c \pm f } n$,
\item $\displaystyle  |S_c(u, v,j)| \in \of{s_c \pm  (c+1)^{-1}f }n^{1/2}$, 
\end{itemize}
 where  $a \pm b$ denotes the interval $[a-b, a+b]$, and the functions $y$,$f$,$q_{b,c}$,$r_c$, and  $s_{c}$ are evaluated at the point $t(j)$.
\end{enumerate}
Note that if the event $\mc{E}_i$ fails (no matter if it fails due to condition \eqref{nobigcod} or \eqref{dynconc}), then $\mc{E}_{i'}$ also fails for all $i' > i$. 
We now show that the first condition  of the event $\mc{E}_i$ hold w.h.p. for every $i$ under consideration. We use the asymptotic equivalence of the models $G(n, m)$ and $G(n, p)$ (where $p = m/\binom{n}{2}$) and the fact that \eqref{nobigcod} is a monotone graph property (see \cite{JLR}). Now to see that this holds w.h.p. we calculate the expected number of pairs $u, v$ with at least $c_{max}$ common neighbors. At step~$i$ the number of edges we have added is at most $n^{3/2} (\log \log n) / 10000$.  Thus it is enough to show that \eqref{nobigcod} holds w.h.p. in $G(n,p)$ where $p\le n^{-1/2} (\log \log n)/5000$.
Now, the expected number of pairs  of vertices in $G(n, p)$ with codegree at least  $c_{max}$ is at most
\begin{align*}
n^2 \binom{n}{c_{max}} p^{2c_{max}} \le n^2 \rbrac{\frac{enp^2}{c_{max}}}^{c_{max}} 
&\le n^2 \rbrac{\frac{(\log \log n)^3}{\log n}}^{c_{max}}
= e^{2\log n} \rbrac{\frac{(\log \log n)^3}{\log n}}^{c_{max}}\\ 
&\le e^{2\log n} \rbrac{\frac{1}{(\log n)^{5/6}}}^{c_{max}} 
= e^{-(\log n)/2} = o(1).
\end{align*}

%


In Sections~\ref{subsec:deg}-\ref{subsec:S} we prove that \eqref{dynconc} also w.h.p. holds.

Since unmatched codegrees are so important to this process, in our analysis we will frequently need to know, for a given pair of vertices $u, v,$ how many possible choices for the next edge $e_i$ would increase $\codeg_U(u, v)$. We denote by $A(u, v)= A(u, v, i)$ the set of such possibilities for $e_i$, which we will now estimate with the assumption that the good event holds.
Suppose 
$w$ is a neighbor of $u$ (resp.~$v$). If we add the edge $\lbrace v, w \rbrace$ (resp. $\lbrace u, w \rbrace$), it  may actually be removed in the same step since it might create a triangle. So as long as we ignore the $\tilde{O}(1)$ vertices in $\codeg(u,v)$, the number of $w$ such that $vw$ is not removed is $S_0(u,v)$ (resp.~$S_0(v,u)$). Thus we have
$$
A(u,v):=S_0(v,u)+S_0(u,v) - \tilde{O}(1).
$$
So for $\alpha(t):= 2s_0$ and all $j \le i$ we have
$$
A(u,v) \in (\alpha \pm 3f)n^{1/2}
$$
for evaluation at $t(j)$. We use $3f$ as the error function here so we can ignore the $\tilde{O}(1)$ term from the definition of $A(u,v)$ (recall that $f = \Omega(n^{-1/5})$). This also allows us to ignore any $\tilde{O}(1)$ edges that might be in $M$ already. Thus we may say $\alpha \le 2$.

In our analysis we will also need to estimate, for a given unmatched edge $e = \lbrace u, v \rbrace \in U(i)$,  a ``count" related to the possibility that $e$ becomes matched in the next step. We put quotes around ``count" because actually we need a weighted count: for each possibility for $e_i=\lbrace v, w\rbrace$ (or $\lbrace u, w \rbrace$) that might result in $e$ becoming matched, we weight it by the probability $1/(c+1)$ that the triangle selected to go into $M(i+1)$ is the triangle containing $e$, where $c$ is the unmatched codegree of $w$ and $u$ (resp. $v$), excluding $v$ (resp. $u$). To this end, we define $K(u,v)$ to be the random variable 
\[
K(u,v) := A(u,v) + \sum_{c = 1}^{c_{max}} \frac{1}{c+1} (S_c(u,v) + S_c(v,u)). 
\]
Here we want to avoid counting the $\tilde{O}(1)$ edges that might be in $M$ already. This, however, is accounted for in the $\tilde{O}(1)$ term in $A(u,v)$. In addition, notice that for $e = \lbrace u, v \rbrace$, the probability that $e$ becomes matched in step $i$ is $$K(u,v)\cdot \frac{2}{n^2} (1~+~\tilde{O}(n^{-1/2})).$$

Now define 
\[
\kappa(y) :=2e^{-y^2}\sum_{c = 0}^{\infty} \frac{y^{2c+1}}{(c+1)!} =  \left\{\begin{array}{lr}
       0, & y=0\\
       2y^{-1}(1-e^{-y^2}), & \text{ otherwise}.
        \end{array}\right.
\]
It is easy to check that $\kappa$ is continuous and twice differentiable for $y \ge 0$. In addition, for $y > 0$,
\begin{align*}
\kappa(y) = 2y^{-1}(1-e^{-y^2}) \le 2y^{-1}(1 - (1-y^2))=2y;
\end{align*}
therefore when $0 \le y \le \zeta$, we have
\begin{equation}
\kappa(y) \le 2\zeta  \le 2.\label{eq:kappa-bnd}
\end{equation}
Then from the dynamic concentration of $S(u,v)$ and $\sum_{c= 0}^{\infty} (c+1)^{-2} =\pi^2/6 < 2$, we have 
\begin{align}
K(u,v)& = \sum_{c = 0}^{c_{max}} \frac{1}{c+1} (S_c(u,v) + S_c(v,u)) - \tilde{O}(1) \nn\\
&\le  \sum_{c = 0}^{c_{max}}\Bigg( \frac{2e^{-y^2}y^{2c+1}}{(c+1)!} + \frac{2f}{(c+1)^{2}} \Bigg)n^{1/2} - \tO{1}\nn\\
&= n^{1/2} \sum_{c = 0}^{\infty} \frac{2e^{-y^2}y^{2c+1}}{(c+1)!} + n^{1/2} \sum_{c = 0}^{c_{max}} \frac{2f}{(c+1)^{2}} + \tO{1}\label{eq:Kuv}\\
&\le \k(y)n^{1/2} + 4f n^{1/2}.\nn
\end{align}

In \eqref{eq:Kuv} we use the fact that for $c \ge c_{max}$ we have $c! \ge \exp\lbrace(3+o(1))\log n\rbrace$ and hence 
$$
n^{1/2}\sum_{c = c_{max}}^{\infty} \frac{y^{2c+1}}{(c+1)!} < n^{-5/2+o(1)} = O(n^{-2}).
$$
A similar argument using the lower bound for $S_c(u,v)$ gives us a lower bound for $K(u,v)$.
Thus
$$
K(x,y) \in (\kappa\pm 4f)n^{1/2},
$$
with evaluation at $t(j)$.

Then straightforward, but somewhat tedious, calculations show that the above functions satisfy the  following differential equations, where $q_{b,c}'$, $r_c'$ and $s_{c}'$ denote 
derivatives of $q_{b,c}$, $r_c$ and $s_{c}$  as functions of~$t$: 
\begin{align} 
q_{b, c}' &= 2q_{b-1,c}\alpha + 2q_{b,c-1}\alpha + 4(b+1)\kappa q_{b+1,c} + 4(c+1)\kappa q_{b,c+1} - 4q_{b,c}(\alpha+b\kappa + c\kappa),\label{eq:qdiff}\\
r_c' &= 2 r_{c-1}\alpha + 4(c+1)\kappa r_{c+1}-(2\alpha+4c\kappa)r_c,\label{eq:rdiff}\\
s_{c}' &= 2s_{c-1}\alpha+4(c+1)\kappa s_{c+1} + 2q_{c,0}-2(\alpha+2c\kappa+\kappa)s_c.\label{eq:sdiff}
\end{align}
These differential equations can be viewed as idealized one-step changes in the random variables $Q_{b,c}(u,v)$, $R_c(v)$, and $S_c(u,v)$. Each of these variables counts copies of some type of substructure, and these copies can be created or destroyed by the process when we add or remove edges. Equations \eqref{eq:qdiff}--\eqref{eq:sdiff} can be understood as expressing the one-step changes in the random variables in terms of these creations and deletions, on average. We will ultimately use these differential equations to argue that the random variables stay close to their deterministic counterparts.

\subsection{Tracking  $d_U(v,j)$}\label{subsec:deg}

First observe that Chernoff's bound implies that w.h.p.
$$ d_{G}(v,j) \in \of{2t \pm n^{-1/4} \log^2 n}n^{1/2}.$$  
Moreover, in order to  estimate  $d_U(v,j)$ it suffices to track  $d_M(v,j)$.

We define the natural filtration $\mc{F}_i$ to be the history of the process up to step $i$. In particular, conditioning on $\mc{F}_i$ tells us the current state of the process. Assuming we are in the event $\mc{E}_{i-1}$, we calculate the expected one-step change of the matched degree, conditional  on 
 $\mc{F}_{i-1}$, namely, 
\[
\Mean{\D d_M(v,i) | \mc{F}_{i-1}} = \Mean{d_M(v,i)-d_M(v,i-1) | \mc{F}_{i-1}}.
\]

We have already revealed $i-1$ edges. Now we  reveal a new edge $e_{i}$.
Note that $d_M(v)$ is nondecreasing. If $e_{i} \subseteq N_U(v)$, where $N_U(v)$ is the set of vertices connected to $v$ in the graph~$U$, then $d_M(v)$ could increase by 2 or not increase at all. For a fixed vertex $u$ in $N_U(v)$, if the edge $e_i = \lbrace u, w\rbrace$ for some $w \in N_U(v)$ with $\codeg_U(u, w) = c+1$, then $d_M(v)$ increases by~2 with probability $1/(c+1)$. The number of such $w$ can be counted with $S_c(u,v)$. Notice that, by taking a sum over $u\in N_U(v)$, we double count such $w$. Finally, if  $e_{i}$ is the edge $\lbrace v, u \rbrace$ for some vertex $u$ not in $N_U(v)$ such that $\codeg_U(u, v) >0$, then $d_M(v)$ increases by $2$. Hence, we have

\begin{align*}
&\Mean{\D d_M (v,i)| \mc{F}_{i-1}}\\
&\qquad = \sbrac{\sum_{u\in N_U(v)}\sum_{c=0}^{c_{max}}2\cdot\frac12(c+1)^{-1} \cdot S_c(u,v, i-1) + \sum_{c=1}^{c_{max}} 2 \cdot R_c(v, i-1)}  \frac {2}{n^2} (1+\tilde{O}(n^{-1/2}))\\
&\qquad \le \Bigg[2(y+f)\sum_{c=0}^{c_{max}} (c+1)^{-1}(s_c + (c+1)^{-1}f) + 4\sum_{c=1}^{c_{max}} (r_c + f)\Bigg]n^{-1} (1+\tilde{O}(n^{-1/2}))\\
&\qquad =\Bigg[2\sum_{c=0}^{c_{max}} (c+1)^{-1}ys_c + 4 \sum_{c = 1}^{c_{max}} r_c + \left(-4 + \sum_{c=0}^{c_{max}} \left(\frac{2y}{(c+1)^2} + \frac{2s_c}{c+1} + 4\right)\right)f \\
&\qquad\qquad \qquad  + f^2 \cdot \sum_{c=0}^{c_{max}} \frac{2}{(c+1)^2}
\Bigg]n^{-1} (1+\tilde{O}(n^{-1/2}))
\end{align*}
 where the functions $y$ and $f$ are evaluated at point $t(i-1)$.
 Now,
\begin{align*}
\sum_{c=0}^{c_{max}} \frac{1}{c+1}y \left(\frac{e^{-y^2}y^{2c+1}}{c!}\right) = e^{-y^2}\sum_{c=0}^{\infty} \frac{y^{2(c+1)}}{(c+1)!} + O(n^{-2}) &= 1-e^{-y^2}+ O(n^{-2})
\end{align*}
where the first equality uses the fact that for $c \ge c_{max}$ we have
\[
c! = \exp\cbrac{(1 + o(1))c \log c} \ge \exp\cbrac{(3+o(1)) \log n}, 
\]
 and so 
\[
 \sum_{c=c_{max}}^\infty \frac{y^{2(c+1)}}{(c+1)!}   < n^{-3+o(1)} = O(n^{-2}). 
\]
In a similar manner,
\begin{align*}
\sum_{c=1}^{c_{max}}  \frac{e^{-y^2}y^{2c}}{c!} = e^{-y^2}\left(\sum_{c=0}^{\infty} \frac{y^{2c}}{c!} - 1\right) + O(n^{-2}) &= 1-e^{-y^2}+ O(n^{-2}).
\end{align*}
Since $y \le 1$ and $s_c \le 1/c!$, we may estimate the coefficient of $f$ as follows:
\begin{align}
-4 +& \sum_{c=0}^{c_{max}} \left(\frac{2y}{(c+1)^2} + \frac{2s_c}{c+1}  + 4\right)  \le -4 + 2\sum_{c=0}^\infty \frac{1}{(c+1)^2} + \sum_{c=0}^\infty \frac{2}{(c+1)!} + 4c_{max} \nn\\
&= -4 +2 \cdot \frac{~\pi^2}{6} + 2(e-1) + 4c_{max} \le -4 + 4 + 4 + 4c_{max} = 4 + 4c_{max}.\nn 
\end{align}
Further, since $2(1-e^{-y^2}) = y\kappa(y) \le 2$ from \eqref{eq:kappa-bnd} and $f^2 \le n^{-1/10}f$, we may write
\begin{align}
&\Mean{\D d_M (v,i)| \mc{F}_{i-1}} \le \Bigg[6-6e^{-y^2} + 4c_{max}f+O(f)\Bigg]n^{-1} +\tilde{O}(n^{-3/2}).\label{eq:Dd}
\end{align}
Define variables 
\[
D^\pm(v)=D^\pm(v, i):=\begin{cases} 
& d_M(v, i) - (2t(i)-y(t(i))  \pm f(t(i)))n^{1/2} \;\;\; \mbox{ if $\mc{E}_{i-1}$ holds}\\
& D^\pm (v, i-1) \;\;\; \mbox{\hskip0.31\textwidth\relax otherwise}.
\end{cases}
\]
We will show that the sequence $D^+(v)$ is a supermartingale. Symmetric calculations show that the $D^-(v)$ is a submartingale.
To do that, we first apply Taylor's theorem to approximate the change in the deterministic function by its derivative. Let $g(t) := 2t-y(t) + f(t)$ and $t(i) := \frac{i}{n^{3/2}}$. Then,
\[
(g\circ t)(i) - (g\circ t)(i-1) = (g\circ t)'(i-1)+\frac{(g\circ t)''(\omega)}{2} = g'(t(i-1))n^{-3/2} + \frac{(g\circ t)''(\omega)}{2},
\]
where $\omega\in [i-1,i]$. But 
\[
(g\circ t)''(i) = ( g'(t(i)) n^{-3/2} )' = g''(t(i)) n^{-3} = (-y''(t)+f''(t))n^{-3}.
\]
Furthermore, by~\eqref{eq:z2} we get that~$|y''(t)|\le 24$.
Also,
\[
f''(t) = \rbrac{\frac{1000\log n}{\log \log n}}^2 \exp\cbrac{\frac{1000\log n}{\log \log n} \cdot t} n^{-1/5} = \rbrac{\frac{1000\log n}{\log \log n}}^2 f(t).
\]
Thus, $(g\circ t)''(\omega) = O(n^{-2})$ and
\begin{equation}\label{eq:g}
(g\circ t)(i) - (g\circ t)(i-1) =  (2-y'(t(i-1)) + f'(t(i-1)))n^{-3/2} + O(n^{-2}).
\end{equation}
Now if we are in $\mc{E}_{i-1}$, then \eqref{eq:Dd} and  \eqref{eq:g} for $t=t(i-1)$ imply
\begin{align*}
\Mean{\D D^+ (v,i)| \mc{F}_{i-1}} &\le \rbrac{ -f' +4c_{max}f +O(f)}n^{-1}  + \tilde{O}(n^{-3/2}) \\
&= \sbrac{-\frac{1000\log n}{\log \log n}+\frac{12\log n}{\log \log n}+O(1) }f n^{-1} + \tilde{O}(n^{-3/2}) \le 0,
\end{align*}
 showing that the sequence $D^+ (v,i)$ is a supermartingale.

We now show that the probability of $D^+ (v)$ becoming positive is small,  implying that there is only a small probability that $d_M(v)$ strays from our desired bounds. 
To do this, we apply now the following martingale inequality due to Freedman:

\begin{lemma}[Freedman~\cite{freedman}] \label{lem:Freedman}
Let $Y(i)$ be a supermartingale with $\Delta Y(i) \leq C$ for all $i$, and let $V(i) :=\displaystyle \sum_{k \le i} \Var[ \Delta Y(k)| \mathcal{F}_{k}]$.  Then,
\[
\P\left[\exists i: V(i) \le b, Y(i) - Y(0) \geq \lambda \right] \leq \displaystyle \exp\left(-\frac{\lambda^2}{2(b+C\lambda) }\right).
\] 
\end{lemma}
 
Observe that $|\D d_M(v, i)| = O(1)$, since at most two edges adjacent to $v$ can become matched at step $i$. Moreover, due to~\eqref{eq:g},
$|\D (2t(i)-y(t(i)) + f(t(i)))n^{1/2}| = O(n^{-1})$ trivially. The triangle inequality thus implies that $|\D D^+ (v,i)| = O(1)$.  Also since the variable $d_M(v, i)$ is nondecreasing we have $\E[|\Delta d_M(v, i)|  | \mathcal{F}_{i}]= \E[\Delta d_M(v, i)  | \mathcal{F}_{i}] = O(n^{-1})$ by \eqref{eq:Dd}.
Hence, the triangle inequality yields $\E[|\Delta D^+|| \mathcal{F}_{k}] = O(n^{-1})$. So the one-step variance is
\[
\Var[ \Delta D^+| \mathcal{F}_{k}] \le \E[(\Delta D^+)^2| \mathcal{F}_{k}] \le O(1) \cdot \E[|\Delta D^+|| \mathcal{F}_{k}] = O(n^{-1}).
\]
Therefore, for Freedman's inequality we 
 use $b = O(n^{-1}) \cdot O(n^{3/2} \log \log n) = \tilde{O}(n^{1/2})$. The ``bad" event here is the event that we have $D^+ (v, i)>0$, and since $D^+(v, 0)=-n^{3/10}$ we set $\lambda=n^{3/10}$. Then, applying Lemma~\ref{lem:Freedman} with $\lambda = n^{3/10}$, $b= \tilde{O}(n^{1/2})$ and $C=\tilde{O}(1)$ yields that the failure probability is at most
\[
\exp \left\{- \frac{n^{3/5}}{ \tilde{O}(n^{1/2}) + \tilde{O}(1) \cdot n^{3/10}}\right\},
\]
which is small enough to beat a union bound over all vertices.

Using symmetric calculations one can apply Freedman's inequality to the supermartingale $-D^-(v, i)$ to show that the ``bad" event $D^- (v, i)<0$ does not occur w.h.p..

\subsection{Tracking $R_c(v)$} 
Since the tracking of $R_c(v)$ seems especially illustrative, we eschew alphabetical order and estimate $\Mean{\D R_c(v,i) | \mc{F}_{i-1}}$ next. Because $R_c(v, i)$ counts the number of vertices $u$  such that $\codeg_U(u,v)=c$, we are interested to know how these codegree functions  can increase or decrease.

\begin{figure}[t]
\begin{subfigure}{.4\textwidth}
\centering
\scalebox{.6}{
\begin{tikzpicture}
[line width = .5pt, 
codeg/.style={ellipse, minimum width=1cm, minimum height=3cm, draw=black!100, fill=gray!0, thick, label=codeg$_U(u{,}v){=}c-1$},
vtx/.style={circle,draw,black,thick,fill=black, line width = .5pt, inner sep=1.5pt},
empty/.style={inner sep=0pt}
]
\FPset{\nbhscale}{1};
\FPset{\nbhy}{4};
\FPset{\nbhx}{0};
\FPset{\xoffset}{.14};
\FPset{\yoffset}{1.4};

\node[empty] (ntv) at (\nbhx+\xoffset,\nbhy+\yoffset) {};
\node[empty] (nbv) at (\nbhx-\xoffset,\nbhy-\yoffset) {};
\node[empty] (ntu) at (\nbhx-\xoffset,\nbhy+\yoffset) {};
\node[empty] (nbu) at (\nbhx+\xoffset,\nbhy-\yoffset) {};
\node[vtx, label=$v$] (v) at (3,0) {};
\node[vtx, label=$u$] (u) at (-3,0) {};
\node[vtx, label=$w$] (w) at (0,1.5) {};

\foreach \b in {ntv, nbv}
    \draw (v) -- (\b);
    
\foreach \b in {ntu, nbu}
    \draw (u) -- (\b);

\draw[] (u) -- (w);
\draw[dotted, thick] (v) -- (w);
\node[codeg] (codeg) at (\nbhx,\nbhy) {};

\end{tikzpicture}
}
\vspace{2ex}
\caption{Ways  $u\in R_{c-1}(v,i-1)$ could join $R_c(v,i)$ with the addition of the dotted edge.}
\label{fig:Rc-a}
\end{subfigure}
\hspace{.1\textwidth}
\begin{subfigure}{.4\textwidth}
\centering
\scalebox{.6}{
\begin{tikzpicture}
[line width = .5pt, 
codeg/.style={ellipse, minimum width=1cm, minimum height=3cm, draw=black!100, fill=gray!0, thick, label=codeg$_U(u{,}v){=}c+1$},
vtx/.style={circle,draw,black,thick,fill=black, line width = .5pt, inner sep=1.5pt},
empty/.style={inner sep=0pt}
]

\FPset{\nbhy}{4};
\FPset{\nbhx}{0};
\FPset{\xoffset}{.14};
\FPset{\yoffset}{1.4};

\node[empty] (ntv) at (\nbhx+\xoffset,\nbhy+\yoffset) {};
\node[empty] (nbv) at (\nbhx-\xoffset,\nbhy-\yoffset) {};
\node[empty] (ntu) at (\nbhx-\xoffset,\nbhy+\yoffset) {};
\node[empty] (nbu) at (\nbhx+\xoffset,\nbhy-\yoffset) {};
\node[vtx, label=$v$] (v) at (3,0) {};
\node[vtx, label=$u$] (u) at (-3,0) {};

\foreach \b in {ntv, nbv}
    \draw (v) -- (\b);
    
\foreach \b in {ntu, nbu}
    \draw (u) -- (\b);

\node[codeg] (codeg) at (\nbhx,\nbhy) {};
\node[vtx, label=$w$] (w) at (0,4) {};
\draw[] (u) -- (w);
\draw[ dashdotted, very thick] (v) -- (w);

\end{tikzpicture}}
\vspace{2ex}
\caption{Ways  $u\in R_{c+1}(v,i-1)$ could join $R_c(v,i)$ with the removal of the dashed edge.}
\label{fig:Rc-b}
\end{subfigure}
\caption{Cases considered in the expected one-step change of $R_{c}(v, i)$.}
\end{figure}
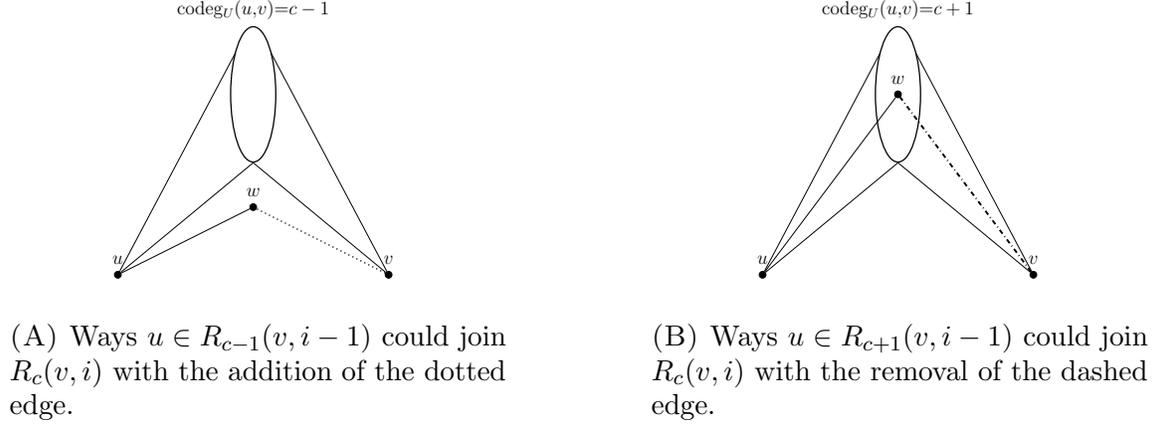

We first focus on the cases that occur the most often. Note that  $\codeg_U(u, v)$ can increase by at most one at any step. Hence, we can increase $\codeg_U(u, v)$ by one if $e_i = \lbrace x, y\rbrace$ such that $x = u$ (resp. $x = v$), $y$ is adjacent to $v$ (resp. $u$), and $e_i$ does not create a triangle with other edges in $U$ (see Figure~\ref{fig:Rc-a}). In the  event $\mc{E}_i$, the number of such edges $e_i$ is $A(u, v)$. Thus, for the positive contribution, we would take the sum of $A(u, v)$ over all $u$ in $R_{c-1}(v)$. But $R_{c}(v)$ can also decrease if $\codeg_U(u,v)$ increases in this way for $u\in R_{c}(v)$. Thus the negative contribution can be described by taking the sum over $A(u, v)$ for all $u$ in $R_{c}(v)$.

The unmatched codegree can also decrease. This can occur when an edge $\lbrace u, w \rbrace$ or $\lbrace v, w \rbrace$ becomes matched, for some $w$ in the unmatched codegree of $u$ and $v$ (see Figure~\ref{fig:Rc-b}). Recall from our previous discussion that $K(w,v)$ and $K(w,u)$ can be used to track these changes. So we obtain a positive contribution by taking the sum of $K(w,v) + K(w,u)$ over all $u$ in $R_{c+1}(v)$ and all $w$ in the unmatched codegree of $u$ and $v$. If $u$ comes from $R_c(v)$, though, the contribution to the unmatched codegree would be negative. 

There are two cases affecting $R_c(v,i-1)$ that occur less frequently. The first such way is to create a triangle by having $e_i$ be in the common neighborhood of $u$ and $v$. The second way is for the edge $e_i$ to contain $u$ and $v$. Either way when we sum over all $u$, the total change is $\tilde{O}(n^{-1})$.

Thus, we obtain 
 \begin{align}
&\Mean{\D R_c(v,i)| \mc{F}_{i-1}}\nn\\
=   &\Bigg[ \sum_{u \in R_{c-1}(v)}  A(u,v) + \sum_{\substack{u \in R_{c+1}(v)\\ w \in \codeg_U(u,v)}} \sbrac{K(w,v) + K(w,u)}   \nn \\
&\quad -\sum_{u \in R_{c}(v)}  A(u,v) - \sum_{\substack{u \in R_{c}(v)\\ w \in \codeg_U(u,v)}} \sbrac{K(w,v) + K(w,u)} \Bigg]\cdot\frac{2}{n^2} (1+ \tilde{O}(n^{-1/2})) + \tilde{O}(n^{-1}). \label{eq:DeltaR}
\end{align}
Estimating this gives
\begin{align}
&\Mean{\D R_c(v,i)| \mc{F}_{i-1}}\nn\\
&\le \bigg[(r_{c-1}+f)(\alpha + 3f)  + 2 (r_{c+1} + f)(c+1)(\kappa + 4f)\nn\\
&\qquad-(r_{c}-f)(\alpha - 3f)  - 2 (r_{c} - f)\cdot c\cdot(\kappa - 4f) \bigg]2n^{3/2}\cdot n^{-2} + \tilde{O}(n^{3/2}\cdot n^{-5/2})+ \tilde{O}(n^{-1})\nn\\ 
&= \bigg[2 r_{c-1}\alpha -(2\alpha+4c\kappa)r_c+ 4(c+1)\kappa r_{c+1} 
 \label{eq:Rdiff}  \\ 
&+8c\kappa f  + \bigg(6r_{c-1}+2(8c+3)r_c + 16(c+1)r_{c+1}+4\alpha + 4\kappa \bigg)f 
+16f^2 
 \bigg]n^{-1/2} + \tilde{O}(n^{-1})\nn
\end{align}
where all functions are evaluated at point $t(i-1)$.

Observe that $8c\kappa f \le 20 c f$ and from the bounds on $r_c, \alpha$ and $\kappa$ in  \eqref{eq:cpqt-bnd} and \eqref{eq:kappa-bnd}, we get that all other terms with $f$ are $O(f)$. Further the $f^2$ terms are also $O\left(f\right)$ since $f^2 \le n^{-1/10} f$. Lastly observe that \eqref{eq:Rdiff} is $r_c'$ by \eqref{eq:rdiff}. 
Thus, $\Mean{\D R_c(v,i)| \mc{F}_{i-1}}$ is at most
\begin{align}
&\bigg[r_c'+ 20 c f +  O(f)\bigg] n^{-1/2} + \tilde{O}(n^{-1}).\label{eq:DC}
\end{align}

Now we define variables
\[
R_c^\pm(v)=R_c^\pm(v, i):=\begin{cases} 
& R_c(v, i) - (r_c(t(i))  \pm f(t(i)))n \;\;\; \mbox{ if $\mc{E}_{i-1}$ holds}\\
& R_c^\pm (v, i-1) \;\;\; \mbox{ \hskip0.29\textwidth\relax otherwise}.
\end{cases}
\]

We show that $R^+(v)$ are supermartingales.
By Taylor's theorem for $g(t) := r_c(t) + f(t)$ and $t(i) := \frac{i}{n^{3/2}}$ we have 
\[
(g\circ t)(i) - (g\circ t)(i-1) = (g\circ t)'(i-1)+\frac{(g\circ t)''(\omega)}{2} = g'(t(i-1))n^{-3/2} + \frac{(g\circ t)''(\omega)}{2},
\]
where $\omega\in [i-1,i]$. But 
\[
(g\circ t)''(i) = ( g'(t(i)) n^{-3/2} )' = g''(t(i)) n^{-3} = r_c''(t)n^{-3}+f''(t)n^{-3}.
\]
However, this is $O(n^{-2})$. To see this, note that 
$$
r_c''(t) = -\frac{(2y^4-(4c+1)y^2+2c^2-c)y^{2c-2}(8e^{-y^2}-12e^{-2y^2})}{c!},
$$
and so $|r_c''(t(i-1))| = O(1)$ along with $|f''(t)| = O(n^{-2})$.
Hence,
$$
\D (r_c(t(i)) + f(t(i)))n 
= \sbrac{r_c'(t(i-1))+f'(t(i-1))} n^{-1/2} +  O(n^{-2}).
$$
Therefore for $t=t(i-1)$ due to \eqref{eq:DC} we get
\begin{align*}
\Mean{\D R_c^+ (v,i)| \mc{F}_{i-1}} 
&\le \Bigg[-f'+20c f+O(f)\Bigg]n^{-1/2} + \tilde{O}(n^{-1})\nn\\
&\le \Bigg[ -\frac{1000 \log n}{\log \log n}+ \frac{60 \log n}{\log \log n} +O(1)\Bigg]fn^{-1/2} + \tilde{O}(n^{-1})\nn \le 0.
\end{align*}

Now observe that $|\D R_c (v)| = \tilde{O}(n^{1/2})$. Indeed, if the new edge $e_i$ has one  vertex at $v$ and the other at say $x$, then this only affects the codegree of $v$ with the $\tilde{O}(n^{1/2})$ many neighbors of $x$. On the other hand if $e_i$ is not incident with $v$ then $v$ loses at most two unmatched edges, say $\lbrace v, x \rbrace$ and $\lbrace v, y \rbrace$, in which case only the codegree of $v$ with the $\tilde{O}(n^{1/2})$ neighbors of $x$ and $y$ can be affected. Thus, we also have $|\D R_c^+ (v)| = \tilde{O}(n^{1/2})$, since $r_c$ and $f$ have much smaller one-step changes. Now we would like to bound $\E[|\Delta R_c(v)|| \mathcal{F}_{k}]$, so we will re-examine~\eqref{eq:DeltaR}. There are positive and negative contributions to $\Delta R_c(v)$, and of course \eqref{eq:DeltaR} represents the expected positive contributions minus the expected negative contributions. Now by the triangle inequality $|\Delta R_c(v)|$ is at most the sum of the positive and negative contributions, and so 
\begin{align}
\E[|\Delta R_c(v)|| \mathcal{F}_{k}]
&\quad \le \Bigg[ \sum_{u \in R_{c-1}(v)}  A(u,v) + \sum_{\substack{u \in R_{c+1}(v)\\ w \in \codeg_U(u,v)}} \sbrac{K(w,v) + K(w,u)}\nn    \\
&\qquad +\sum_{u \in R_{c}(v)}  A(u,v) + \sum_{\substack{u \in R_{c}(v)\\ w \in \codeg_U(u,v)}} \sbrac{K(w,v) + K(w,u)} \Bigg]\cdot\frac {2}{n^2} + \tilde{O}(n^{-1}) 
\nn\\
&\quad=O(n^{-1/2}), \label{eq:Dabs}
\end{align}
since each term in \eqref{eq:Rdiff} is $O(n^{-1/2})$ as demonstrated in~\eqref{eq:DC}.
Thus,
\[
\E[|\Delta R_c^+(v)|| \mathcal{F}_{k}] \le \E[|\Delta R_c(v)|| \mathcal{F}_{k}] + |\D (r_c(t)  + f(t))|n = O(n^{-1/2}),
\]
and hence the one-step variance is 
\[
\Var[ \Delta R_c^+(v)| \mathcal{F}_{k}] \le \E[(\Delta R_c^+(v))^2| \mathcal{F}_{k}] = \tilde{O}(n^{1/2}) \cdot \E[|\Delta R_c^+(v)|| \mathcal{F}_{k}] = \tilde{O}(1).
\]

The ``bad" event here is the event that $R_c^+ (v, i)>0$. Since $R_c^+(v, 0)=n^{4/5}$ we set $\lambda=n^{4/5}$. Then, Lemma~\ref{lem:Freedman} applied with $\lambda = n^{4/5}$,  $b= \tilde{O}(n^{3/2})$ and $C=\tilde{O}(n^{1/2})$ yields that the failure probability is at most
\[
\exp \left\{- \frac{\tilde{O}(n^{8/5})}{ \tilde{O}(n^{3/2}) + \tilde{O}(n^{1/2})\cdot n^{4/5}}\right\},
\]
which is small enough to beat a union bound over all vertices as well as possible values of $c$.

\subsection{Tracking $Q_{b, c}(u,v)$}\label{subsec:Q}

\begin{figure}[t]
\begin{subfigure}{.4\textwidth}
\centering
\scalebox{.6}{
\begin{tikzpicture}
[line width = .5pt, 
codeg/.style={ellipse, minimum width=1cm, minimum height=3cm, draw=black!100, fill=gray!0, thick, label=codeg$_U(u{,}w) {=}b-1$},
codegc/.style={ellipse, minimum width=1cm, minimum height=3cm, draw=black!100, fill=gray!0, thick, label=below:codeg$_U(v{,}w){=}c$},
vtx/.style={circle,draw,black,thick,fill=black, line width = .5pt, inner sep=1.5pt},
empty/.style={inner sep=0pt}
]
\FPset{\nbhy}{4};
\FPset{\nbhx}{0};
\FPset{\xoffset}{.14};
\FPset{\yoffset}{1.4};
\FPset{\nbhyc}{-5};
\FPset{\nbhxc}{0};
\FPset{\xoffsetc}{.14};
\FPset{\yoffsetc}{1.4};

\node[empty] (ntv) at (\nbhx+\xoffset,\nbhy+\yoffset) {};
\node[empty] (nbv) at (\nbhx-\xoffset,\nbhy-\yoffset) {};
\node[empty] (ntu) at (\nbhx-\xoffset,\nbhy+\yoffset) {};
\node[empty] (nbu) at (\nbhx+\xoffset,\nbhy-\yoffset) {};
\node[empty] (ntvc) at (\nbhxc+\xoffset,\nbhyc+\yoffset) {};
\node[empty] (nbvc) at (\nbhxc-\xoffset,\nbhyc-\yoffset) {};
\node[empty] (ntuc) at (\nbhxc-\xoffset,\nbhyc+\yoffset) {};
\node[empty] (nbuc) at (\nbhxc+\xoffset,\nbhyc-\yoffset) {};
\node[vtx, label=$w$] (w) at (3,-.5) {};
\node[vtx, label=$u$] (u) at (-3,0) {};
\node[vtx, label=$v$] (v) at (-3,-1) {};

\foreach \b in {ntv, nbv}
    \draw (w) -- (\b);
    
\foreach \b in {ntu, nbu}
    \draw (u) -- (\b);
    
\foreach \b in {ntvc, nbvc}
    \draw (v) -- (\b);
    
\foreach \b in {ntuc, nbuc}
    \draw (w) -- (\b);    

\node[vtx, label=$w'$] (wp) at (0,1.5) {};
\draw[] (u) -- (wp);
\draw[dotted, thick] (w) -- (wp);
\node[codeg] (codeg) at (\nbhx,\nbhy) {};
\node[codegc] (codegc) at (\nbhxc,\nbhyc) {};

\end{tikzpicture}
}
\vspace{2ex}
\caption{Ways $w\in Q_{b-1, c}(u,v,i-1)$ could join $Q_{b, c}(u,v,i)$ with the addition of the dotted edge.}
\label{fig:Qbc-a}
\end{subfigure}
\hspace{.1\textwidth}
\begin{subfigure}{.4\textwidth}
\centering
\scalebox{.6}{
\begin{tikzpicture}
[line width = .5pt,
codeg/.style={ellipse, minimum width=1cm, minimum height=3cm, draw=black!100, fill=gray!0, thick, label=codeg$_U(u{,}w) {=}b+1$},
codegc/.style={ellipse, minimum width=1cm, minimum height=3cm, draw=black!100, fill=gray!0, thick, label=below:codeg$_U(v{,}w){=}c$},
vtx/.style={circle,draw,black,thick,fill=black, line width = .5pt, inner sep=1.5pt},
empty/.style={inner sep=0pt}
]
\FPset{\nbhy}{4};
\FPset{\nbhx}{0};
\FPset{\xoffset}{.14};
\FPset{\yoffset}{1.4};
\FPset{\nbhyc}{-5};
\FPset{\nbhxc}{0};
\FPset{\xoffsetc}{.14};
\FPset{\yoffsetc}{1.4};

\node[empty] (ntv) at (\nbhx+\xoffset,\nbhy+\yoffset) {};
\node[empty] (nbv) at (\nbhx-\xoffset,\nbhy-\yoffset) {};
\node[empty] (ntu) at (\nbhx-\xoffset,\nbhy+\yoffset) {};
\node[empty] (nbu) at (\nbhx+\xoffset,\nbhy-\yoffset) {};
\node[empty] (ntvc) at (\nbhxc+\xoffset,\nbhyc+\yoffset) {};
\node[empty] (nbvc) at (\nbhxc-\xoffset,\nbhyc-\yoffset) {};
\node[empty] (ntuc) at (\nbhxc-\xoffset,\nbhyc+\yoffset) {};
\node[empty] (nbuc) at (\nbhxc+\xoffset,\nbhyc-\yoffset) {};
\node[vtx, label=$w$] (w) at (3,-.5) {};
\node[vtx, label=$u$] (u) at (-3,0) {};
\node[vtx, label=$v$] (v) at (-3,-1) {};

\foreach \b in {ntv, nbv}
    \draw (w) -- (\b);
    
\foreach \b in {ntu, nbu}
    \draw (u) -- (\b);
    
\foreach \b in {ntvc, nbvc}
    \draw (v) -- (\b);
    
\foreach \b in {ntuc, nbuc}
    \draw (w) -- (\b);    

\node[codeg] (codeg) at (\nbhx,\nbhy) {};
\node[codegc] (codegc) at (\nbhxc,\nbhyc) {};
\node[vtx, label=$w'$] (wp) at (0,4) {};
\draw[] (u) -- (wp);
\draw[ dashdotted, very thick] (w) -- (wp);

\end{tikzpicture}
}
\vspace{2ex}
\caption{Ways $w\in Q_{b+1, c}(u,v,i-1)$ could join $Q_{b, c}(u,v,i)$ with the removal of the dashed edge.}
\label{fig:Qbc-b}
\end{subfigure}
\caption{Cases considered in the expected one-step change of $Q_{b,c}(u,v, i)$.}
\end{figure}

We continue by calculating $\Mean{\D Q_{b, c}(u,v,i)| \mc{F}_{i-1}}$. As before, we consider positive and negative contributions, beginning with the positive ones. First, for $w\in Q_{b-1, c}(u,v)$, $e_i$ could increase $\codeg_U(u,w)$ up to $b$ from $b-1$ (see Figure~\ref{fig:Qbc-a}). This situation is accounted with $A(u,w)$. Second, we could have a symmetric case for $w\in Q_{b,c-1}(u,v)$ tracked by  $A(v,w)$. Third, for $w\in Q_{b+1, c}(u,v)$ an edge $\lbrace u, w' \rbrace$ (or $\lbrace w, w' \rbrace$) could become matched for $w'\in\codeg_U(u,w)$. This would result in $\codeg_U(u,w) = b$ (see Figure~\ref{fig:Qbc-b}). This can be tracked with $K(u,w') + K(w,w')$. Fourth, a similar case to the third can occur but with $w\in Q_{b,c+1}(u,v)$. 

For the negative contribution, all cases are symmetric. The only difference is that all $w$ are taken from $Q_{b,c}(u,v)$ and thus these gains or losses of edges destroy the existing structure.

There are also three unlikely cases that affect $Q_{b, c}(u,v)$. In each case, we consider $w\in Q_{b',c'}(u,v)$ for with $b'$ and $c'$ appropriately chosen from $ \lbrace b-1, b, b+1\rbrace$ and $\lbrace c-1, c, c+1\rbrace$, respectively. There are $\tilde{O}(n)$ such $w$. First, the edge $e_i$ could be in the common neighborhood of $u, w$ or of $v, w$. Second, $e_i$ could be either the edge $\lbrace v, w \rbrace$ or $\lbrace u, w \rbrace$, reducing the codegree of the vertices in that edge. Third, $e_i$ could have one vertex in the common neighborhood of $u, w$ and the other in the common neighborhood of $v, w$. In all of these cases, we have a change of $\tilde{O}(1)$ and hence together with the edge probability, we have a change of $\tilde{O}(n^{-1})$.

Consequently,
\begin{align}
&\!\!\!\!\!\Mean{\D Q_{b, c}(u,v,i)| \mc{F}_{i-1}} \nn\\
&=    \Bigg[ \sum_{w \in Q_{b-1, c}(u, v)} A(u, w) +  \sum_{w \in Q_{b, c-1}(u, v)} A(v, w)  +\sum_{\substack{w \in Q_{b+1, c}(u, v) \\ w'\in \codeg(u,w)}} (K(u, w')+K(w,w')) \nn\\
&\qquad+\sum_{\substack{w \in Q_{b, c+1}(u, v) \\ w'\in \codeg(v,w)}} (K(v, w')+K(w,w')) -\sum_{w \in Q_{b, c}(u, v)} A(u, w) - \sum_{w \in Q_{b, c}(u, v)} A(v, w)  \nn\\
&\qquad-\sum_{\substack{w \in Q_{b, c}(u, v) \\ w'\in \codeg(u,w)}} (K(u, w')+K(w,w'))  -\sum_{\substack{w \in Q_{b, c}(u, v) \\ w'\in \codeg(v,w)}} (K(v, w')+K(w,w')) \Bigg]\nn\\
&\qquad\qquad\qquad\times \frac {2}{n^2} (1+\tilde{O}(n^{-1/2})) +\ \tilde{O}(n^{-1}) \label{eq:DeltaQ}\\
&\le \Bigg[(q_{b-1,c}+ f)(\alpha+3f) + (q_{b,c-1}+ f)(\alpha+3f) 
+ 2(\kappa+4f)(q_{b+1,c}+ f)(b+1) \nn\\
&\qquad+ 2(\kappa+4f)(q_{b,c+1}+ f)(c+1)
-2(q_{b,c}- f)(\alpha-3f)    \nn \\
& \qquad-2(b+c)(\kappa-4f)(q_{b,c}- f) 
\Bigg]\cdot 2n^{-1/2} + \tilde{O}(n^{-1}) \nn
\end{align}
\vspace{-.8cm}
\begin{align}
&=\Bigg[2q_{b-1,c}\alpha + 2q_{b,c-1}\alpha + 4q_{b+1,c}(b+1)\kappa + 4q_{b,c+1}(c+1)\kappa - 4q_{b,c}(\alpha+b\kappa + c\kappa)  \label{eq:Qdiff}\\
&\qquad+ 8\kappa (b+c) f+\bigg(6q_{b-1,c}+6q_{b,c-1}+16(b+1)q_{b+1,c}+16(c+1)q_{b,c+1}\nn\\
&\qquad\qquad+16q_{b,c}(b+ c) + 12q_{b,c}+8(\alpha + \kappa) \bigg)f
+32f^2\Bigg]
\cdot n^{-1/2} + \tilde{O}(n^{-1}).\nn
\end{align}
We simplify this. First notice that $ 8\kappa (b+c) f \le 20 (b+c)f$ and
from the bounds $q_{b,c} \le 1/(b!c!)$ and $\alpha, \kappa \le 3$, we see that all other multiples of $f$ are $O(f)$.
Also recall that $f^2\le n^{-1/10}f = O(f)$. Finally, observe that by \eqref{eq:qdiff}, line \eqref{eq:Qdiff} is $q_{b,c}'$.
Therefore, $\Mean{\D Q_{b, c}(u,v,i)| \mc{F}_{i-1}}$ is at most
\begin{align}
&\bigg[ q_{b,c}'+20(b+c)f + O(f)\bigg] n^{-1/2} + \tilde{O}(n^{-1}).\label{eq:DQ}
\end{align}
Now we define variables
\[
Q_{b, c}^\pm(u,v)=Q_{b, c}^\pm(u,v, i):=\begin{cases} 
& Q_{b, c}(u,v, i) - (q_{b, c}(t(i))  \pm f(t(i)))n \;\;\; \mbox{ if $\mc{E}_{i-1}$ holds}\\
& Q_{b, c}^\pm (u,v, i-1) \;\;\; \mbox{ otherwise}.
\end{cases}
\]

We demonstrate that $Q_{b,c}^+(u,v)$ is a supermartingale by using Taylor's theorem applied to $g(t) := q_{b,c}(t) +  f(t)$ and $t(i) := \frac{i}{n^{3/2}}$. In this case,
\[
(g\circ t)(i) - (g\circ t)(i-1) = (g\circ t)'(i-1)+\frac{(g\circ t)''(\omega)}{2} = g'(t(i-1))n^{-3/2} + \frac{(g\circ t)''(\omega)}{2},
\]
where $\omega\in [i-1,i]$. But 
\[
(g\circ t)''(i) = ( g'(t(i)) n^{-3/2} )' = g''(t(i)) n^{-3} = q_{b,c}''(t)n^{-3}+f''(t)n^{-3}.
\]
However,
$$
q_{b,c}''(t) = -\frac{16\left(e^{y^2}-\frac32\right)\left(4y^4-4(b+c+1)y^2+(b+c)\left(b+c-\frac12\right)\right)e^{-3y^2}y^{2b+2c-2}}{b!c!},
$$
so $|q_{b,c}''(t(i-1))| = O(1)$. Also, we have  $| f''(t)| = O(n^{-2})$.
Hence, 
\begin{align*}
&\D (q_{b,c}(t(i)) + f(t(i)))n
= \sbrac{q_{b,c}'(t(i-1))+f'(t(i-1))} n^{-1/2} +  O(n^{-2}).
\end{align*}
Therefore from \eqref{eq:DQ} with $t=t(i-1)$, we have
\begin{align*}
\Mean{\D Q_{b,c}^+ (u,v,i)| \mc{F}_{i-1}} 
&\le \Bigg[-f'+20(b+c)f+O(f)\Bigg]n^{-1/2} + \tilde{O}(n^{-1})\nn\\
&\le \Bigg[ -\frac{1000 \log n}{\log \log n}+ \frac{120\log n}{\log \log n} +O(1)\Bigg]fn^{-1/2} + \tilde{O}(n^{-1})\nn \le 0.
\end{align*}

Let us consider the effect on $Q_{b, c} (u,v)$ by removing one edge $e$ from the unmatched graph. If $e$ is incident with $u$, say $e=\lbrace u, x\rbrace$, then the only vertices $w \in Q_{b, c} (u,v)$ that could be affected are in the set $\lbrace x\rbrace \cup N(x)$ which has size $\tilde{O}(n^{1/2})$. Similarly if $e$ is incident with $v$. If $e$ is not incident with $u, v$ then the only affected $w \in Q_{b, c} (u,v)$ would be the endpoints of $e$.  Thus we have $|\D Q_{b, c} (u,v)| = \tilde{O}(n^{1/2})$,  and also $|\D Q_{b, c}^+ (u,v)| = \tilde{O}(n^{1/2})$  because the deterministic terms in $ Q_{b, c}^+(u, v)$ have much smaller one-step changes.  We can also see that $\E[|\Delta Q_{b, c}(u,v)|| \mathcal{F}_{k}]= O(n^{-1/2})$ by another argument analogous to the one used to justify~\eqref{eq:Dabs}. Indeed, $\E[|\Delta Q_{b, c}(u,v)|| \mathcal{F}_{k}]$ is at most the sum of the absolute values of the terms in \eqref{eq:DeltaQ}, all of which are $O(n^{-1/2})$. Thus,
\[
\E[|\Delta Q_{b, c}^+(u,v)|| \mathcal{F}_{k}] \le \E[|\Delta Q_{b, c}(u,v)|| \mathcal{F}_{k}] + |\D (q_{b,c}(t)  + f(t))n| = O(n^{-1/2}), 
\]
 and the one-step variance is 
\[
\Var[ \Delta Q_{b, c}^+(u, v)| \mathcal{F}_{k}] \le \E[(\Delta Q_{b, c}^+(u, v))^2| \mathcal{F}_{k}] = \tilde{O}(n^{1/2}) \cdot \E[|\Delta Q_{b, c}^+(u, v)|| \mathcal{F}_{k}] = \tilde{O}(1).
\]
Thus, Lemma \ref{lem:Freedman} applied with $\lambda=n^{4/5}$, $\tilde{b}=\tilde{O}(n^{3/2})$ and $C=\tilde{O}(n^{1/2})$ (using $\tilde{b}$ instead of $b$ to avoid notational collision) yields that  the failure probability is at most
\[
\exp \left\{- \frac{n^{8/5}}{  \tilde{O}(n^{3/2})  +  \tilde{O}(n^{1/2}\cdot n^{4/5}) }\right\} 
\]
which is  again small enough to beat a union bound over all pairs of vertices and values of~$b,c$.

\subsection{Tracking $S_c(u,v)$}\label{subsec:S}

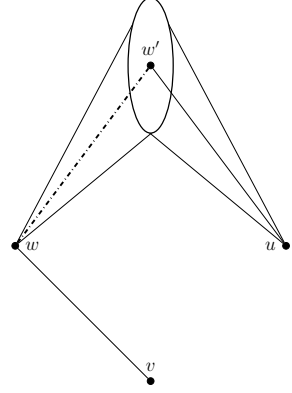
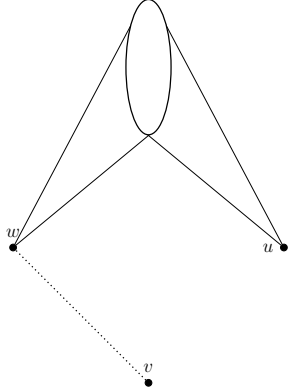
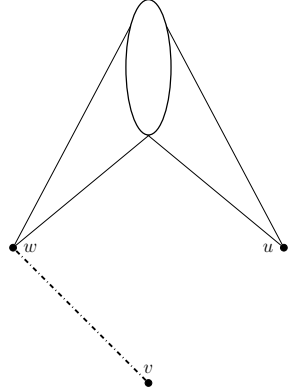
\begin{figure}[h]
\begin{subfigure}{.4\textwidth}
\centering
\scalebox{.6}{
\begin{tikzpicture}
[line width = .5pt, 
codeg/.style={ellipse, minimum width=1cm, minimum height=3cm, draw=black!100, fill=gray!0, thick, label=codeg$_U^*(u{,}w){=}c-1$},
vtx/.style={circle,draw,black,thick,fill=black, line width = .5pt, inner sep=1.5pt},
empty/.style={inner sep=0pt},
nbh/.style={ellipse, minimum width=7cm, minimum height=1cm, draw=black!100, fill=gray!0, thin}
]
\FPset{\nbhscale}{1};
\FPset{\nbhy}{4};
\FPset{\nbhx}{0};
\FPset{\y}{0};
\FPset{\x}{0};
\FPset{\xoffset}{.14};
\FPset{\yoffset}{1.4};
\FPset{\xoffsetv}{3.5};
\FPset{\yoffsetv}{.04};

\node[empty] (ntv) at (\nbhx+\xoffset,\nbhy+\yoffset) {};
\node[empty] (nbv) at (\nbhx-\xoffset,\nbhy-\yoffset) {};
\node[empty] (ntu) at (\nbhx-\xoffset,\nbhy+\yoffset) {};
\node[empty] (nbu) at (\nbhx+\xoffset,\nbhy-\yoffset) {};

\node[vtx, label=left:$u$] (u) at (3,0) {};
\node[vtx, label=right:$w$] (w) at (-3,0) {};
\node[vtx, label=$w'$] (wp) at (0,1.5) {};

\foreach \b in {ntv, nbv}
    \draw (u) -- (\b);
    
\foreach \b in {ntu, nbu}
    \draw (w) -- (\b);

\draw[] (u) -- (wp);
\draw[dotted, thick] (wp) -- (w);
\node[vtx, label=$v$] (v) at (0, -3){};

\foreach \b in {w}
    \draw (v) -- (\b);

\node[codeg] (codeg) at (\nbhx,\nbhy) {};

\end{tikzpicture}
}
\vspace{2ex}
\caption{Ways  $w\in S_{c-1}(u,v,i-1)$ could join $S_{c}(u,v,i)$ with the addition of the dotted edge.}
\label{fig:Sc-a}
\end{subfigure}
\hspace{.1\textwidth}
\begin{subfigure}{.4\textwidth}
\centering
\scalebox{.6}{
\begin{tikzpicture}
[line width = .5pt, 
codeg/.style={ellipse, minimum width=1cm, minimum height=3cm, draw=black!100, fill=gray!0, thick, label=codeg$_U^*(u{,}w){=}c+1$},
vtx/.style={circle,draw,black,thick,fill=black, line width = .5pt, inner sep=1.5pt},
empty/.style={inner sep=0pt},
nbh/.style={ellipse, minimum width=7cm, minimum height=1cm, draw=black!100, fill=gray!0, thin}
]
\FPset{\nbhscale}{1};
\FPset{\nbhy}{4};
\FPset{\nbhx}{0};
\FPset{\y}{0};
\FPset{\x}{0};
\FPset{\xoffset}{.14};
\FPset{\yoffset}{1.4};
\FPset{\xoffsetv}{3.5};
\FPset{\yoffsetv}{.04};

\node[empty] (ntv) at (\nbhx+\xoffset,\nbhy+\yoffset) {};
\node[empty] (nbv) at (\nbhx-\xoffset,\nbhy-\yoffset) {};
\node[empty] (ntu) at (\nbhx-\xoffset,\nbhy+\yoffset) {};
\node[empty] (nbu) at (\nbhx+\xoffset,\nbhy-\yoffset) {};

\node[vtx, label=left:$u$] (u) at (3,0) {};
\node[vtx, label=right:$w$] (w) at (-3,0) {};

\foreach \b in {ntv, nbv}
    \draw (u) -- (\b);
    
\foreach \b in {ntu, nbu}
    \draw (w) -- (\b);

\foreach \b in {w}
    \draw (v) -- (\b);

\node[codeg] (codeg) at (\nbhx,\nbhy) {};
\node[vtx, label=$w'$] (wp) at (0,4) {};
\draw[] (u) -- (wp);
\draw[ dashdotted, very thick] (w) -- (wp);
\node[vtx, label=$v$] (w) at (0, -3){};

\end{tikzpicture}
}
\vspace{2ex}
\caption{Ways  $w\in S_{c+1}(u,v,i-1)$ could join $S_{c}(u,v,i)$ with the removal of the dashed edge.}
\label{fig:Sc-b}
\end{subfigure}
\begin{subfigure}{.4\textwidth}
\vspace{2ex}
\centering
\scalebox{.6}{

\begin{tikzpicture}
[line width = .5pt, 
codeg/.style={ellipse, minimum width=1cm, minimum height=3cm, draw=black!100, fill=gray!0, thick, label=codeg$_U^*(u{,}w){=}c+1$},
vtx/.style={circle,draw,black,thick,fill=black, line width = .5pt, inner sep=1.5pt},
empty/.style={inner sep=0pt},
nbh/.style={ellipse, minimum width=7cm, minimum height=1cm, draw=black!100, fill=gray!0, thin}
]
\FPset{\nbhscale}{1};
\FPset{\nbhy}{4};
\FPset{\nbhx}{0};
\FPset{\y}{0};
\FPset{\x}{0};
\FPset{\xoffset}{.14};
\FPset{\yoffset}{1.4};
\FPset{\xoffsetv}{3.5};
\FPset{\yoffsetv}{.04};

\node[empty] (ntv) at (\nbhx+\xoffset,\nbhy+\yoffset) {};
\node[empty] (nbv) at (\nbhx-\xoffset,\nbhy-\yoffset) {};
\node[empty] (ntu) at (\nbhx-\xoffset,\nbhy+\yoffset) {};
\node[empty] (nbu) at (\nbhx+\xoffset,\nbhy-\yoffset) {};

\node[vtx, label=left:$u$] (u) at (3,0) {};
\node[vtx, label=$w$] (w) at (-3,0) {};

\foreach \b in {ntv, nbv}
    \draw (u) -- (\b);
    
\foreach \b in {ntu, nbu}
    \draw (w) -- (\b);

\node[codeg] (codeg) at (\nbhx,\nbhy) {};
\draw[dotted, thick] (v) -- (w);
\node[vtx, label=$v$] (w) at (0, -3){};

\end{tikzpicture}
}
\vspace{2ex}
\caption{Ways  $w\in Q_{c,0}(u,v,i-1)$ could join $S_{c}(u,v,i)$.}
\label{fig:Sc-c}
\end{subfigure}
\hspace{.1\textwidth}
\begin{subfigure}{.4\textwidth}
\vspace{2ex}
\centering

\scalebox{.6}{

\begin{tikzpicture}
[line width = .5pt, 
codeg/.style={ellipse, minimum width=1cm, minimum height=3cm, draw=black!100, fill=gray!0, thick, label=codeg$_U^*(u{,}w){=}c$},
vtx/.style={circle,draw,black,thick,fill=black, line width = .5pt, inner sep=1.5pt},
empty/.style={inner sep=0pt},
nbh/.style={ellipse, minimum width=7cm, minimum height=1cm, draw=black!100, fill=gray!0, thin}
]
\FPset{\nbhscale}{1};
\FPset{\nbhy}{4};
\FPset{\nbhx}{0};
\FPset{\y}{0};
\FPset{\x}{0};
\FPset{\xoffset}{.14};
\FPset{\yoffset}{1.4};
\FPset{\xoffsetv}{3.5};
\FPset{\yoffsetv}{.04};

\node[empty] (ntv) at (\nbhx+\xoffset,\nbhy+\yoffset) {};
\node[empty] (nbv) at (\nbhx-\xoffset,\nbhy-\yoffset) {};
\node[empty] (ntu) at (\nbhx-\xoffset,\nbhy+\yoffset) {};
\node[empty] (nbu) at (\nbhx+\xoffset,\nbhy-\yoffset) {};

\node[vtx, label=left:$u$] (u) at (3,0) {};
\node[vtx, label=right:$w$] (w) at (-3,0) {};

\foreach \b in {ntv, nbv}
    \draw (u) -- (\b);
    
\foreach \b in {ntu, nbu}
    \draw (w) -- (\b);

\node[empty] (wl) at (\x+\xoffsetv,\y-\yoffsetv) {};
\node[empty] (ur) at (\x-\xoffsetv,\y-\yoffsetv) {};

\node[codeg] (codeg) at (\nbhx,\nbhy) {};

\draw[dashdotted, very thick] (v) -- (w);

\node[vtx, label=$v$] (w) at (0, -3){};

\end{tikzpicture}
}
\vspace{2ex}
\caption{Ways to lose $w\in S_{c}(u,v,i-1)$.}
\label{fig:Sc-d}
\end{subfigure}
\vspace{2ex}
\caption{Cases considered in the expected one-step change of $S_c(u,v, i)$. Note that codeg$_U^*(u{,}w)$ denotes the unmatched codegree of $u$ and $w$, excluding~$v$.}
\end{figure}

Finally, we calculate $\Mean{\D S_{c}(u,v,i)| \mc{F}_{i-1}}$, starting with positive contributions. Denote by $\codeg_U^*(u,w)$ the unmatched codegree of $u$ and $w$, excluding $v$. First, for $w\in S_{c-1}(u,v)$, we could increase $\codeg_U^*(u,w)$ by the addition of an edge, a situation accounted for by $A(u,w)$ (see Figure~\ref{fig:Sc-a}). Second, for $w\in S_{c+1}(u,v)$, an edge $\lbrace w, w' \rbrace$ or $\lbrace u, w' \rbrace$ could become matched, thereby reducing $\codeg_U^*(u,w)$ from $c+1$ to $c$. This situation is handled with $K(w,w')+K(u,w')$ and shown in Figure~\ref{fig:Sc-b}. Third, $e_i$ could be an edge $\lbrace v, w \rbrace$ for some $w$ with $\codeg_U^*(w,u)=c$. We may use $Q_{c,0}(u,v)$ for this case (see Figure~\ref{fig:Sc-c}).

Next, we consider negative contributions. The first two cases are identical to the positive contribution with the exception that we consider $w\in S_c(u,v)$, and hence lose such a $w$ instead of gaining one. For the third case, the edge $\lbrace w, v \rbrace$ could become matched for $w\in S_c(u,v)$ and thus $w$ no longer meets the requirements for membership in $S_c(u,v)$. This is shown in Figure~\ref{fig:Sc-d}. 

We are also concerned with more unlikely cases. Take $w\in S_{c'}(u,v)$ for appropriate $c'$. There are $\tilde{O}(n^{1/2})$ of such $w$. First, the chosen edge could be in the common neighborhood of $u,w$. Second, $e_i$ could be the edge $\lbrace u, w \rbrace$, which would form triangles with the common neighborhood of these vertices. Both of these are a change of $\tilde{O}(1)$, and so the over all change for these unlikely cases is $\tilde{O}(n^{-3/2})$.

Therefore,
\begin{align}
&\!\!\!\Mean{\D S_{c}(u,v,i)| \mc{F}_{i-1}} \nn\\
&=  \Bigg[ \sum_{w\in S_{c-1}(u,v)} A(u,w) + Q_{c,0}(u,v) + \sum_{\substack{w\in S_{c+1}(u,v) \\ w'\in \codeg_U(u,w)}} \sbrac{K(w,w') + K(u, w')}-\sum_{w\in S_c(u,v)} A(u,w)\nn \\
&  - \sum_{\substack{w\in S_{c}(u,v) \\ w'\in \codeg_U(u,w)}} \sbrac{K(w,w') + K(u, w')} - \sum_{w\in S_c(u,v)} K(v,w) \Bigg] \cdot \frac {2}{n^2} (1+\tilde{O}(n^{-1/2})) +\ \tilde{O}(n^{-3/2}) \label{eq:DeltaS} \\
&\le \bigg[(s_{c-1}+c^{-1}f)(\alpha+3f) + (q_{c,0} + f) + 2(c+1)(s_{c+1}+(c+2)^{-1}f)(\kappa + 4f)\nn\\
&\hspace{5cm}-(s_{c}-(c+1)^{-1}f)(\alpha-3f)-2c(s_c-(c+1)^{-1}f)(\kappa-4f)\nn\\
&\hspace{5cm}-(s_c-(c+1)^{-1}f)(\kappa-4f)\bigg]\cdot 2n^{-1} +\tilde{O}(n^{-3/2}) \nn\\
&\le \bigg[2s_{c-1}\alpha+4(c+1)\kappa s_{c+1} + 2q_{c,0}-2(\alpha+2c\kappa+\kappa)s_c\label{eq:Sdiff}\\
&\quad+\big(2 + 4c\kappa(c+2)^{-1} + 4c\kappa(c+1)^{-1}\big)f \nn\\
&\quad+\big(2\alpha c^{-1} + 16(c+1) s_{c+1} + 2(8c+7) s_c + 6s_{c-1}    
 + 4\kappa(c+2)^{-1} + 2(\alpha + \kappa) (c+1)^{-1}\big)f\label{eq:SRerr1}\\
&\quad+\big(16(c+1)(c+2)^{-1}-16c(c+1)^{-1}+6c^{-1} - 14(c+1)^{-1}\big)f^2
\bigg] n^{-1} + \tilde{O}(n^{-3/2}).\nn 
\end{align}
Recall from \eqref{eq:sdiff} that \eqref{eq:Sdiff} is $s_c'$. We also remark that
\begin{align*}
2 + 4c\kappa(c+2)^{-1} + 4c\kappa(c+1)^{-1} \le 50.
\end{align*}
From the bounds on $\alpha$, $s_c$ and $\kappa$ in \eqref{eq:cpqt-bnd} and \eqref{eq:kappa-bnd} we can easily conclude that \eqref{eq:SRerr1} is ${O((c+1)^{-1}f)}$. Furthermore, since $f^2 \le n^{-1/10}f$, we can say that the $f^2$ terms are also $O((c+1)^{-1}f)$. Thus,
\begin{align}
&\Mean{\D S_{c}(u,v,i)| \mc{F}_{i-1}}
\le \big[s_c' + 50f+O((c+1)^{-1}f)\big]n^{-1}+ \tilde{O}(n^{-3/2}). \label{eq:DS}
\end{align}
Now we define variables
\[
S_{c}^\pm(u,v)=S_{c}^\pm(u,v, i):=\begin{cases} 
& S_{c}(u,v, i) - (s_{c}(t(i))  \pm (c+1)^{-1}f(t(i)))n^{1/2} \;\;\; \mbox{ if $\mc{E}_{i-1}$ holds}\\
& S_{c}^\pm (u,v, i-1) \;\;\; \mbox{ otherwise}.
\end{cases}
\]
Applying Taylor's theorem with $g(t) := s_{c}(t) + (c+1)^{-1} f(t)$ and $t(i) := \frac{i}{n^{3/2}}$ yields
\[
(g\circ t)(i) - (g\circ t)(i-1) = (g\circ t)'(i-1)+\frac{(g\circ t)''(\omega)}{2} = g'(t(i-1))n^{-3/2} + \frac{(g\circ t)''(\omega)}{2},
\]
where $\omega\in [i-1,i]$. But 
\[
(g\circ t)''(i) = ( g'(t(i)) n^{-3/2} )' = g''(t(i)) n^{-3} = s_{c}''(t)n^{-3}+(c+1)^{-1}f''(t)n^{-3}.
\]
However,
$$
s_{c}''(t) = -\frac{-4y^{2c-1}e^{-2y^2}(2y^4-(4c+3)y^2+2c^2+c)(2e^{y^2}-3)}{c!},
$$
so $|s_{c}''(t(i-1))| = O(1)$. Also, we have  $|(c+1)^{-1} f''(t)| = O(n^{-2})$.
Hence, 
\begin{align*}
&\D (s_{c}(t(i)) + (c+1)^{-1}f(t(i)))n^{1/2}
= \sbrac{s_{c}'(t(i-1))+(c+1)^{-1}f'(t(i-1))} n^{-1} +  O(n^{-3/2}).
\end{align*}
So as a consequence of \eqref{eq:DS} for $t=t(i-1)$ we have 
\begin{align*}
\Mean{\D S_{c}^+ (u,v,i)| \mc{F}_{i-1}} 
&\le \Bigg[-(c+1)^{-1}f' + 50f+O((c+1)^{-1}f)\Bigg]n^{-1} + \tilde{O}(n^{-3/2})\nn\\
&\le \Bigg[ {-\frac{1000 \log n}{\log \log n}}+50(c+1) +O(1)\Bigg](c+1)^{-1}fn^{-1} + \tilde{O}(n^{-3/2})\nn\\
&\le \Bigg[ -\frac{1000 \log n}{\log \log n}+ \frac{200\log n}{\log \log n} +O(1)\Bigg](c+1)^{-1}fn^{-1} + \tilde{O}(n^{-3/2})\nn \le 0.
\end{align*}

Next we demonstrate that $|\D S_c (u,v)| = {O}(\log n)$. First, an edge $e$ in in the unmatched graph $G_U$  might be removed. Indeed, if $e$ is incident with $u$, say $e=\lbrace u, x\rbrace$, then the removal of $e$ can only affect vertices $w\in S_c(u, v)$ such that $w \in \{x\} \cup (N(x)\cap N(v))$ of which there are only $O(\log n)$. Similarly if $e$ is incident with $v$ then at most $O(1)$ vertices $w\in S_c(u, v)$ are affected. Finally, if $e$ is not incident with $u, v$ then the only vertices $w\in S_c(u, v)$ that could be affected are the endpoints of $e$. Thus, zero or three edges are removed at any step and each one affects $O(\log n)$ vertices $w$. We can also make symmetrical observations for adding an edge to $G_U$. Hence, $|\D S_c (u,v)| = {O}(\log n)$. Also $|\D S_c^+ (u,v)| = O(\log n)$, since the deterministic terms have much smaller one-step changes. We can also see that $\E[|\Delta S_c(u,v)|| \mathcal{F}_{k}]= O(n^{-1})$ by an argument analogous to the one used to justify \eqref{eq:Dabs}. Indeed, $\E[|\Delta S_c(u,v)|| \mathcal{F}_{k}]$ is at most the sum of the absolute values of the terms in \eqref{eq:DeltaS}, all of which are $O(n^{-1})$. Thus, 
\[
\E[|\Delta S_c^+(u,v)|| \mathcal{F}_{k}] \le \E[|\Delta S_c(u,v)|| \mathcal{F}_{k}] + |\D (s_c(t)  + f(t))n^{1/2}| = O(n^{-1})
\]
and
\[
\Var[ \Delta S_c^+(u, v)| \mathcal{F}_{k}] \le \E[(\Delta S_c^+(u, v))^2| \mathcal{F}_{k}] = O(\log n) \cdot \E[|\Delta S_c^+(u, v)|| \mathcal{F}_{k}] = \tilde{O}(n^{-1}).
\]
Therefore, using Lemma \ref{lem:Freedman} with $\lambda=(c+1)^{-1}n^{3/10}$, $b=\tilde{O}(n^{1/2})$ and $C=\tilde{O}(1)$ our failure probability is at most 
\[
\exp \left\{- \frac{(c+1)^{-2}n^{3/5}}{  \tilde{O}(n^{1/2}) +  \tilde{O}(n^{3/10})}\right\},
\]
which is small enough to beat a union bound over all pairs of vertices and values of $c$. 

\subsection{Proof of Theorem~\ref{thm:main}}

At the end of our process, after revealing $kn^{3/2}$ edges, the number of unmatched edges is at most $\frac{n^{3/2}}{2}(y(\C)+f(\C))$, and thus the number of matched edges is at least 
$$
\C n^{3/2} - \frac{n^{3/2}}{2}(y(k)+f) \ge \C n^{3/2} - \frac{y(\C)n^{3/2}}{2}-n^{7/5}.
$$
Recall that the only edges of $M$ are those of edge-disjoint triangles. Hence, the number of edge-disjoint triangles at the end of our process is w.h.p. at least 
$$
(1+o(1))\frac{1}{3}\left(\C - \frac{y(\C)}{2} \right)  n^{3/2}.
$$

\section{Proof of Theorem \ref{thm:tuza}}\label{sec:tuza}

We bound $\tau(G)$ using two different approaches for different ranges of $\C$. First, we consider the triangle-free process for $t = O(1)$. This process accepts a set of edges forming a triangle-free subgraph of $G(n, m)$ and so the rejected edges form a triangle cover. We will refer to Bohman's original triangle-free paper \cite{bohman}.  Recall that in this process one maintains a triangle-free subgraph $G_{T}(i) \subseteq G(n, i)$  by revealing one edge at a time, and adding that edge to $G_{T}(i)$ only if it does not create a triangle in $G_{T}(i)$. 

In order to use the results about the triangle-free process, we must reconcile the differences between our step parameter and the one used by Bohman. In \cite{bohman}, one step was counted as the acceptance of a proposed edge whereas we count each step as the presenting of an edge. To avoid notational confusion with Bohman's paper, all of the variables in \cite{bohman} will appear in this paper with a circumflex above them, e.g. $\hat{i}, \hatt, \hat{Q},$ etc.
With this notation, the number of edges accepted by the process after $i=tn^{3/2}$ edges are proposed is $\hati = \hatt n^{3/2}$. 

Bohman proved that w.h.p.~for all $\hati \le \C  n^{3/2}$ the number $\hat{Q}(\hati)$ of edges eligible to be inserted into the triangle-free graph (i.e. edges that would be accepted if proposed) is 
\begin{equation}\label{eqn:bohman}
   \hat{Q}(\hati) \in (1\pm n^{-\beta})\binom{n}{2}e^{-4\hatt^2}. 
\end{equation}
where the constant $\beta > 0$ is derived from Bohman's original error function. Since \eqref{eqn:bohman} holds for some $\beta>0$ we may take $\beta$ to be arbitrarily small. We will assume that $0 < \beta <  1/2$.
We also note that Bohman proved this for all $\hati$ at most some constant times $n^{3/2} \log^{1/2}n$ but we will not fully use that here. 

Let the random variable $A(i)$ be the number of edges that have been accepted after $i$ edges have been presented. Of course Bohman's step number is $\hat{i} = A(i)$ and so $\hatt = A(i)n^{-3/2}$. We will track $A(i)$ using the differential equation method. We will show that $A(i) \approx n^{3/2} a(t)$ for some function $a(t)$. 

Let $\tilde{\mc{E}}_i$ be the event that for all $i^* \le i$ we have:
\begin{enumerate}[(i)]
\item $\hat{Q}(A(i^*)) \in (1\pm n^{-\beta})\binom{n}{2}e^{-4 \rbrac{ A(i^*)n^{-3/2}}^2}$,
\item  $A(i^*) \in (1 \pm f_A(t^*))n^{3/2} a(t^*)$.
\end{enumerate}
Since in this section we have $t = O(1)$, we also have $a(t) = O(1)$. Assume the constant $c$ is an upper bound on $a(t)$. We define the error function $f_A$ to be 
$$
f_A(t):=n^{-\beta}e^{(2+13c^2)t} = O(n^{-\beta}).
$$ 
Since Bohman showed that the first condition holds w.h.p., we verify that the second condition also w.h.p. holds. Note that in the event $\tilde{\mc{E}}_i$ we have 
\begin{align*}
\hat{Q}(A(i)) &\in (1\pm n^{-\beta})\binom{n}{2}e^{-4 (1 \pm 3f_A(t)) a(t)^2 }\\
&\subseteq  (1\pm n^{-\beta}) \cdot ( 1 \pm 12f_A(t) a^2(t) + O(f_A^2(t)a^{4}(t))) \cdot \binom{n}{2}e^{-4 a(t)^2 }\\
& \subseteq  \bigg(1\pm n^{-\beta} \pm 12c^2 f_A(t) +  O(n^{-2\beta})\bigg) \binom{n}{2}e^{-4 a(t)^2 }\\
& \subseteq  \bigg(1\pm (1+12c^2) f_A(t) +  O(n^{-2\beta})\bigg) \binom{n}{2}e^{-4 a(t)^2 }
\end{align*}
where the second inclusion follows from the Taylor series for $e^{ 12f_A(t) a^2(t) }$, and the last inclusion follows from $f_A(t) \ge n^{-\beta}$. 
Therefore we have
\[
\E[\Delta A(i) | \mc{F}_i] = \frac{\hat{Q}(A(i))}{\binom{n}{2} - i}
\]
and we may derive the following differential equation:
$$
\frac{da}{dt} = e^{-4a(t)^2} \quad\text{ and }\quad a(0)=0.
$$
Now estimate  $\E[\Delta A(i) | \mc{F}_i]$ as follows, noting that the lower bound can be attained with symmetric calculations:
\begin{align*}
\E[\Delta A(i) | \mc{F}_i] &= \frac{\hat{Q}(A(i))}{\binom{n}{2} - i} \\
& \le \frac{ (1 + (1+12c^2)f_A(t)) \binom{n}{2}e^{-4 a(t)^2 }}{\binom{n}{2} - i}\\
&\le (1+{O}(n^{-1/2}))(1+ (1+12c^2)f_A(t))e^{-4a^2}\\
&\le e^{-4a^2}\big(1+ (1+12c^2)f_A(t)\big)+ {O}(n^{-1/2}).
\end{align*}
Next we define variables which we will show are supermartingales 
\[
A^\pm=A^\pm(i):=\begin{cases} 
& A(i) - (a(t(i))  \pm f_A(t(i)))n^{3/2} \;\;\; \mbox{ if $\tilde{\mc{E}}_{i-1}$ holds}\\
& A^\pm (i-1) \;\;\; \mbox{ otherwise}.
\end{cases}
\]
The calculation to verify $A^+$ is shown. Applying Taylor's theorem with $g(t) := a(t) + f_A(t)$ and $t(i) := \frac{i}{n^{3/2}}$ yields
\[
(g\circ t)(i) - (g\circ t)(i-1) = (g\circ t)'(i-1)+\frac{(g\circ t)''(\omega)}{2} = g'(t(i-1))n^{-3/2} + \frac{(g\circ t)''(\omega)}{2},
\]
where $\omega\in [i-1,i]$. But 
\[
(g\circ t)''(i) = ( g'(t(i)) n^{-3/2} )' = g''(t(i)) n^{-3} = \left(a''(t)+f_A''(t)\right)n^{-3}
\]
as well as
\[
a''(t) = 2a(t) \cdot a'(t) \cdot e^{-4a^2(t)} = 2a(t) e^{-8a^2(t)} = O(1)
\]
and
\[
f_A''(t) = (2+13c^2)^2n^{-\beta}e^{(2+13c^2)t} = O(n^{-\beta}).
\]
This gives us $(g\circ t)''(i) = O(n^{-3})$, so
\begin{align*}
&\D (a(t(i)) + f_A(t(i)))n^{3/2}
= a'(t(i-1))+f_A'(t(i-1))  +  O(n^{-3/2}).
\end{align*}
Thus for $t=t(i-1)$ we have 
\begin{align*}
\Mean{\D A^+ (i)| \mc{F}_{i-1}} 
&\le (1+12c^2)e^{-4a^2}f_A - f_A' + O(n^{-3/2})\\
& \le (1+12c^2)f_A - f_A' + O(n^{-3/2})\\
&= -(1+c^2)e^{(2+13c^2)t}n^{-\beta} + O(n^{-3/2}) \le 0.
\end{align*}
Now we use the Azuma-Hoeffding inequality to complete this analysis. 
\begin{lemma}[Azuma~\cite{azuma}, Hoeffding~\cite{hoeffding}] \label{lem:AH}
If $Y_0, Y_1,\ldots$ variables are supermartingales and w.h.p. $|Y_j-Y_{j-1}| \le C$, then for all positive integers $m$ and $\lambda$
$$
\P(Y_m - Y_0 \ge \lambda) \le \exp\left(-\frac{\lambda^2}{2C^2m}\right).
$$
\end{lemma}
Clearly $|\D A(i)| = O(1)$  and 
\begin{align*}
&|\D (a(t(i)) + f_A(t(i)))n^{3/2}|
\le |a'(t(i-1))|+|f_A'(t(i-1))|  +  O(n^{-3/2}) = O(1),
\end{align*}
and thus for $t = t(i)$, $\lambda = f_A(0)n^{3/2} = n^{13/10}$ and $C = O(1)$,
\begin{align*}
\P(A(i)>n^{3/2}(a(t)) + f_A(t)) &=\P(A^+(i) > 0)\\
&= \P(A^+(i) - A^+(0) > n^{13/10}) \le \exp\left(-\frac{n^{13/5}}{2C^2tn^{3/2}}\right) = o(1).
\end{align*}
Hence w.h.p. $A(i^*) = (1 +o(1))n^{3/2} a(t^*)$ for $t \le t^*$. So we may say that $\tau(G(n, tn^{3/2})) \le (1+o(1))(t-a(t))n^{3/2}$. 

To optimize our upper bound, we combine what we just obtained with another bound on $\tau(G)$ which is better for certain values of $\C$. We observe that one can cover all triangles in any graph $G$ by using at most half of the edges. To demonstrate this, let $H$ be the largest bipartite subgraph of $G$. It is well-known that $|E(H)| \ge \frac{1}{2}|E(G)|$ (see e.g., \cite{Erdos}). Therefore $E(G)\setminus E(H)$ cover all triangles and we have 
\begin{equation*}
    \tau(G(n, m)) \le m/2.
\end{equation*}
Therefore, we can conclude that w.h.p.
\[
    \tau(G(n, \C n^{3/2})) \le (1+o(1))U_\tau(\C) n^{3/2},
\]
where
\begin{equation} 
U_{\tau}(\C) := \min\lbrace\C-a(\C), \C/2 \rbrace.\label{eq:U}
\end{equation}
And from Theorem~\ref{thm:main}, we have 
$$
\nu(G(n, \C n^{3/2})) \ge (1+o(1)) L_\nu^*(\C) n^{3/2},
$$
where
\begin{equation}
L_\nu^*(\C):=  \frac{1}{3}\left(\C - \frac{y(\C)}{2} \right). \label{eq:Lstar}
\end{equation}
By the proof of Theorem~\ref{thm:tuza-old}, to verify Tuza's conjecture for $G(n, \C n^{3/2})$ it suffices to check $\C$ in the range of  $0.2 \le \C \le  3$. However, our bounds \eqref{eq:U} and \eqref{eq:Lstar} are enough to show that the quotient $U_\tau(\C)/L_\nu^*(\C) \le 2$ for this range (see Appendix for details).

\section{Concluding remarks}\label{sec:remarks}

As a possible area of future study, we conjecture that we can improve Tuza's conjecture for the random graph $G(n,m)$. Recall that in general, Baron and Kahn~\cite{BK} showed that we cannot decrease the multiplicative constant $2$ since for any $\alpha>0$ there are arbitrarily large graphs $G$ of positive density satisfying $\tau(G) > (1-o(1))|G|/2$ and $\nu(G) < (1+\alpha)|G|/4$. However, the graphs used in this are significantly different than $G(n,m)$, being a subgraph of a blowup construction.

Frankl and R\"{o}dl proved in~\cite{FR86}  that w.h.p. for every $\varepsilon_1 > 0$ there is some $\C$ such that  the largest triangle-free subgraph of  $G(n, p)$ for $p = \C n^{-1/2}$  has at least $1/2 - \varepsilon_1$ proportion of the edges. In addition, for $\varepsilon_2 > 0$ the triangle packing number $\nu(G(n,m))$ with $m = \C n^{3/2}$ and $\C \ge (\log n)^2$ is w.h.p. $1/3(1-\varepsilon_2)\C n^{3/2}$.  Then by the asymptotic equivalence of $G(n,p)$ and $G(n,m)$ we see that for $G= G(n,m)$, we have
$$
\frac{\tau(G)}{\nu(G)} \le \frac{(1/2 - \varepsilon_1)\C n^{3/2}}{1/3(1 - \varepsilon_2)\C n^{3/2}} \approx \frac32
$$
On the other hand we have $\tau(G) \ge \nu(G)$ for all $G$ and the ratio $\tau(G)/\nu(G)$ is about 1 when triangles start emerging in $G$. We believe this ratio starts at about $1$ and grows to about $3/2$. This leads us to the following conjecture.

\begin{conjecture}
For all $C> 3/2$ and $G = G(n,m)$, w.h.p. $\tau(G) \le C \cdot \nu(G)$ for all range of $m$.
\end{conjecture}


\appendix
\section{Calculations} \label{sec:sysdiff}

Here we verify that $U_\tau(\C)/L_\nu^*(\C) \le 2$. First define $g(k) := (k/2) / L_\nu^*(k)$  and $h(k) := (k-a(k))/L_\nu^*(k)$. We show that 
\begin{enumerate}[(i)]
\item $h$ is increasing for $.2 \le \C \le 1.29$,\label{itm:hinc}
\item $g$ is decreasing for $1.28 \le \C \le 3$, and \label{itm:gdec}
\item $h(1.29), g(1.28) < 2$.\label{itm:hg2}
\end{enumerate}
We first demonstrate \eqref{itm:gdec} and \eqref{itm:hg2} for $g$. Recall that $y(k) \le 0.6368$ and hence $$L_\nu^*(1.28) \ge \frac{1}{3}\left(1.28 - \frac{0.6368}{2} \right) \ge 0.3205.$$ Thus
$$
g(1.28) = \frac{(1.28/2) }{ L_\nu^*(1.28)} \le \frac{0.64}{0.3205} = 1.9969,
$$
verifying \eqref{itm:hg2} for $g$. To show \eqref{itm:gdec}, using $y' = 6e^{-y^2}-4$ gives us
$$
\frac{dg}{dk} = \frac{18 k e^{-y(k)^2} - 3y(k) - 12k}{(-2k + y(k))^2}.
$$
Then the numerator of $dg/dk$ is $0$ when $k = 0$. Now taking the derivative of this numerator tells us that 
$$
18e^{-y^2} - 36k y e^{-y^2} ( 6e^{-y^2}-4) - 3 ( 6e^{-y^2}-4) - 12 = - 36k y e^{-y^2} ( 6e^{-y^2}-4) \le 0
$$
for $\C \ge 0$ since $y', y \ge 0$. Thus, $dg/dk$ is nonpositive for $1.28 \le \C \le 3$.

Verifying \eqref{itm:hinc} and \eqref{itm:hg2} for $h$ analytically, however, seems more complicated.
Therefore, to avoid tedious calculations we use Maple and verify the required conditions numerically.
For $a$, we obtain
\begin{quote}
\texttt{
\hspace{-.45cm} DE\_A := \{diff(a(t), t) = exp(-4*a(t)\textasciicircum 2), a(0) = 0\};\\
a\_sol\_list := dsolve(DE\_A, numeric, output = listprocedure);\\
a\_sol := rhs(a\_sol\_list[2]);
}
\end{quote}
\begin{figure}[h]
\centering
\includegraphics[scale=.6]{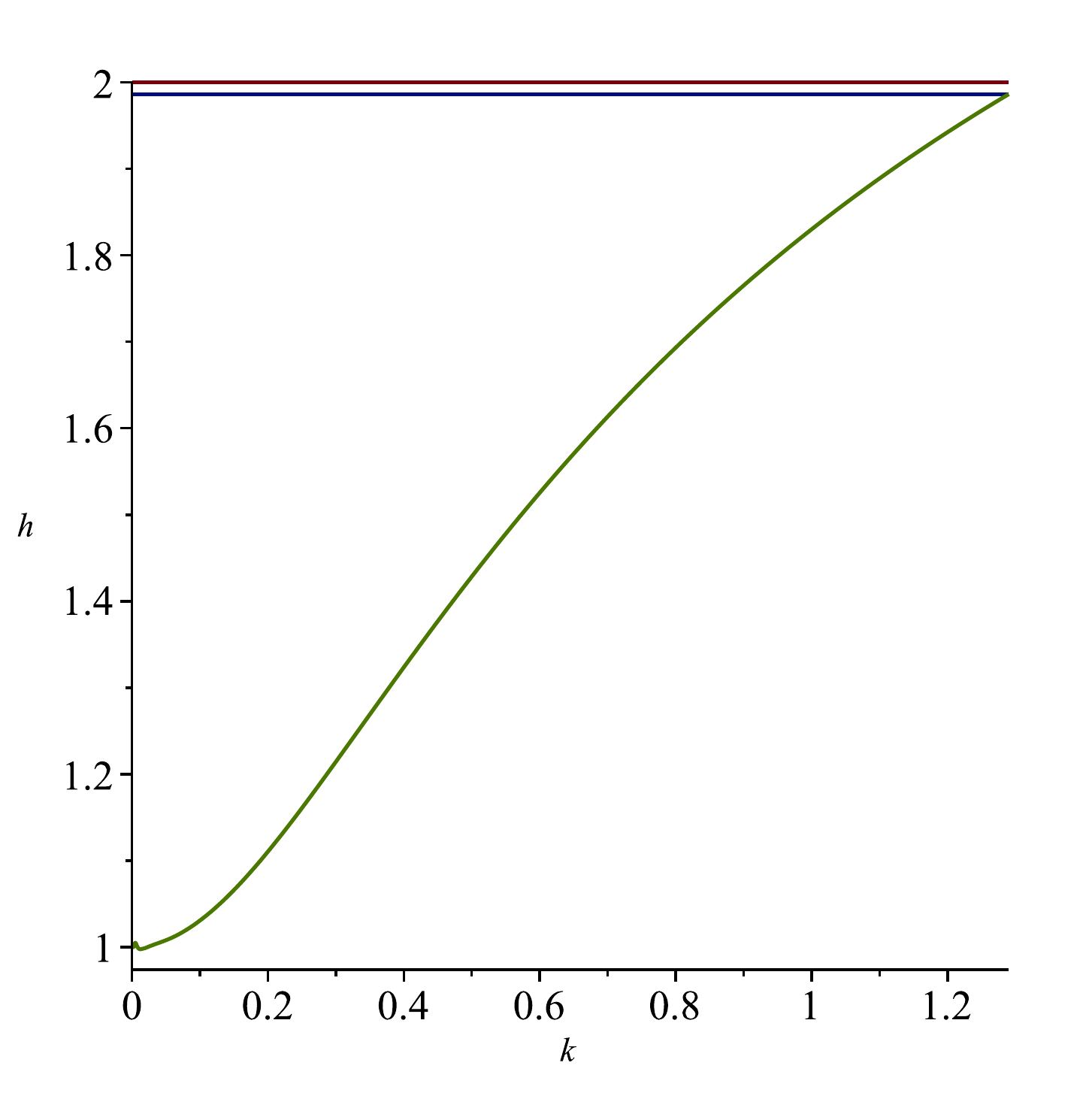}
\vspace{1ex}
\caption{The graph of $h$ (in green) compared with $h(1.29)$ (in blue).}
\label{fig:dh}
\label{fig:hplot}
\end{figure}
With this in hand, we can see that \eqref{itm:hg2} holds since $h(1.29) \le 1.987$. Now plotting the graph of $h$ (see Figure~\ref{fig:hplot}) verifies \eqref{itm:hinc}.


\begin{thebibliography}{99}

\bibitem{AZ} R. Aharoni and S. Zerbib, {\em A generalization of Tuza's conjecture}, to appear in Journal of Graph Theory.

\bibitem{AY} N.~Alon and R.~Yuster,
{\em On a hypergraph matching problem},
Graphs Combin. \textbf{21} (2005), no. 4, 377--384.

\bibitem{azuma} K.~Azuma, {\em Weighted sums of certain dependent random variables}, Tohoku Math. J. \textbf{19} (1967), 357--
367.


\bibitem{BK} J.~Baron and  J.~Kahn,
{\em Tuza's conjecture is asymptotically tight for dense graphs},
Combin. Probab. Comput.~\textbf{25} (2016), no. 5, 645--667.

\bibitem{BG} A.~Basit and D.~Galvin, personal communication.

\bibitem{BD} P.~Bennett and A.~Dudek, {\em A gentle introduction to the differential equation method and dynamic concentration}, submitted.
Available online: \texttt{arXiv:2007.01994}. 

\bibitem{BDZ} P.~Bennett, A.~Dudek, and S.~Zerbib, {\em Large triangle packings and Tuza's conjecture in sparse random graphs}, to appear in Comb. Probab. Comput.  Available online:  \texttt{arXiv:1810.11739}.

\bibitem{bohman} T. Bohman, {\em The triangle-free process,}
 Advances in Mathematics {\bf 221}, (2009) 1653--1677. 

\bibitem{BFL} T. Bohman, A. Frieze, E. Lubetzky, {\em Random triangle removal,}
 Advances in Mathematics {\bf 280} (2015), 379--438.
 
 \bibitem{BK2} T. Bohman and P. Keevash,
{\em Dynamic concentration of the triangle-free process}, Seventh European Conference on Combinatorics, Graph Theory and Applications \textbf{16} (2013),  489--495.

\bibitem{B1} B. Bollob\'as, {\em The life and work of Paul Erd\H{o}s}, Wolf Prize in mathematics. Vol. 1 (S. S. Chern and F. Hirzebrunch, eds.), World Scientific Publishing Co. Inc., River Edge, NJ, 2000, 292--315.

\bibitem{B2} B. Bollob\'as, {\em To prove and conjecture: Paul Erd\H{o}s and his mathematics}, Amer. Math. Monthly {\bf 105} (1998), no. 3, 209--237.

\bibitem{BR} B. Bollob\'as and O. Riordan, {\em Random graphs and branching processes, in: Handbook of large-scale random networks,} Bolyai Soc. Math. Stud. {\bf 18}, Springer, Berlin, 2009, 15--115.


\bibitem{Erdos} P. Erd\H{o}s, On some extremal problems in graph theory, Israel J. Math. 3 (1965) 113--116. 

\bibitem{ESW} P. Erd\H{o}s, S. Suen and P. Winkler, {\em On the size of a random maximal graph}, Random Structures \& Algorithms {\bf 6} (1995), no. 2-3, 309--318. 

\bibitem{FGM} G. Fiz Pontiveros, S. Griffiths and R. Morris,  
{\em The triangle-free process and $R(3,k)$}, Mem. Amer. Math. Soc. \textbf{263} (2020), 125pp. 

\bibitem{FR} P.~Frankl and V.~R\"{o}dl,
{\em Near perfect coverings in graphs and hypergraphs},
European J. Combin. \textbf{6} (1985), no. 4, 317--326.

\bibitem{FR86} P.~Frankl and V.~R\"{o}dl,
{\em Large triangle-free subgraphs in graphs without $K_4$},
Graphs Combin. \textbf{2} (1986), 135--144.


\bibitem{freedman} D.A. Freedman, {\em On Tail Probabilities for Martingales}, Ann. Probability~\textbf{3} (1975), 100--118.

\bibitem{haxell} P.~Haxell,  {\em Packing and covering triangles in graphs}, Discrete Math. {\bf 195} (1999), 251--254.

\bibitem{HR} P.~Haxell and V.~R\"{o}dl, {\em Integer and fractional packings in dense graphs}, Combinatorica~\textbf{21} (2001), 13--38.

\bibitem{hoeffding} W.~Hoeffding, {\em Probability inequalities for sums of bounded random variables}, J. Amer. Statist. Assoc.
\textbf{58} (1963), 13--30.

\bibitem{KP} J.~Kahn and J.~Park, {\em Tuza's conjecture for random graphs}, submitted.

\bibitem{krivelevich} M.~Krivelevich,  {\em On a conjecture of Tuza about packing and covering of triangles},
Discrete Math. \textbf{142}, 281--286.

\bibitem{JLR} S. Janson, T. \L uczak, A. Ruci\'nski, {\em Random Graphs}, Wiley-Interscience, New York (2009).

\bibitem{Makai} T. Makai, {\em The Reverse H-free Process for Strictly 2-Balanced Graphs}, Journal of Graph Theory, {\bf 79} (2015), 125--144.

\bibitem{tuza} Zs. Tuza,  {\em A Conjecture}, Finite and Infinite Sets, Eger, Hungary 1981, A. Hajnal, L. Lov\'asz, V.T. S6s
(Eds.),{\em   Proc. Colloq. Math. Soc. J. Bolyai}, {\bf 37}, North-Holland, Amsterdam, 1984, p. 888.

\bibitem{LutzRev} L. Warnke, {\em On the Method of Typical Bounded Differences}, Comb. Probab. Comput., {\bf 25} (2016), 269--299.

\bibitem{nick2}
N. Wormald, The differential equation method for random graph processes and greedy algorithms, in \emph{Lectures on Approximation and Randomized Algorithms} (M.~Karo\'nski and H.J. Pr\"{o}mel, eds), pp. 73--155. PWN, Warsaw, 1999.

\bibitem{yuster} R.~Yuster, {\em Dense graphs with a large triangle cover have a large triangle packing},
Combin. Probab. Comput. \textbf{21} (2012), 952--962.

\end{thebibliography}
\end{document}